\documentclass[12pt]{amsart}
\usepackage{amscd,amsmath,amsthm,amssymb}
\usepackage[pagewise]{lineno}
\usepackage{amsmath, amsthm, amscd, amsfonts, amssymb, enumerate, verbatim, newlfont, calc, graphicx, color}
\usepackage[bookmarksnumbered, colorlinks=false, plainpages]{hyperref}
\usepackage{doi}
\usepackage{graphicx}
\usepackage[margin=2.5cm]{geometry}
\usepackage{tikz}
\usetikzlibrary{arrows.meta, positioning}
\usepackage{caption}
\usepackage[utf8]{inputenc}
\usepackage{subcaption}
\usepackage{colortbl}
\usepackage{xcolor}
\usepackage{array}
\usepackage{ytableau}

\usepackage{arydshln}

\DeclareUnicodeCharacter{02BC}{'}

\newcolumntype{H}{>{\setbox0=\hbox\bgroup}c<{\egroup}@{}}
\newcolumntype{Z}{>{\setbox0=\hbox\bgroup}c<{\egroup}@{\hspace*{-\tabcolsep}}}

\usepackage[capitalize, noabbrev]{cleveref}
\usepackage{pstricks}
\usepackage{lmodern}
\usepackage{pst-node}
\usepackage[normalem]{ulem}
\psset{unit=0.7cm,linewidth=0.8pt,arrowsize=2.5pt 4}

\newpsstyle{fatline}{linewidth=1.5pt}
\newpsstyle{fyp}{fillstyle=solid,fillcolor=verylight}
\definecolor{verylight}{gray}{0.97}
\definecolor{light}{gray}{0.93}
\definecolor{medium}{gray}{0.82}

\newcommand{\qand}{\quad \mbox{and} \quad}

\newcommand{\qforall}{\quad \mbox{for all} \quad}

\def\NZQ{\mathbb}               
\def\NN{{\NZQ N}}

\def\ZZ{{\NZQ Z}}

\def\MN{{\mathcal N}}

\def\MS{{\mathcal S}}

\def\opn#1#2{\def#1{\operatorname{#2}}} 
\opn\chara{char} \opn\length{\ell} \opn\pd{pd} \opn\rk{rk}
\opn\projdim{proj\,dim} \opn\injdim{inj\,dim} \opn\rank{rank}
\opn\depth{depth} \opn\grade{grade} \opn\height{height}
\opn\embdim{emb\,dim} \opn\codim{codim}

\opn\Tr{Tr} \opn\bigrank{big\,rank}
\opn\superheight{superheight}\opn\lcm{lcm}
\opn\trdeg{tr\,deg}
	\opn\reg{reg} \opn\lreg{lreg} \opn\ini{in} \opn\lpd{lpd}
	\opn\size{size} \opn\sdepth{sdepth}
	\opn\link{link} \opn\del{del}\opn\fdepth{fdepth}\opn\lex{lex}\opn\dist{dist}
	\opn\div{div} \opn\Div{Div} \opn\cl{cl} \opn\Cl{Cl}
	\opn\reg{reg} \opn\lreg{lreg} \opn\ini{in} \opn\lpd{lpd}
	\opn\size{size} \opn\sdepth{sdepth}
	\opn\link{link}\opn\fdepth{fdepth}\opn\lex{lex}\opn\dist{dist}
	\opn\div{div} \opn\Div{Div} \opn\cl{cl} \opn\Cl{Cl}
	\opn\Spec{Spec} \opn\Supp{Supp} \opn\supp{supp} \opn\Sing{Sing}
	\opn\Ass{Ass} \opn\Min{Min}\opn\Mon{Mon}
	\opn\Ann{Ann} \opn\Rad{Rad} \opn\Soc{Soc}
	\opn\Im{Im} \opn\Ker{Ker} \opn\Coker{Coker} \opn\Am{Am}
	\opn\Hom{Hom} \opn\Tor{Tor} \opn\Ext{Ext} \opn\End{End}
	\opn\Aut{Aut} \opn\id{id}
	
	\opn\nat{nat}
	\opn\pff{pf}
	\opn\Pf{Pf} \opn\GL{GL} \opn\SL{SL} \opn\mod{mod} \opn\ord{ord}
	\opn\Gin{Gin} \opn\Hilb{Hilb}\opn\sort{sort}
	\opn\aff{aff} \opn
	\con{conv} \opn\relint{relint} \opn\st{st}
	\opn\lk{lk} \opn\cn{cn} \opn\core{core} \opn\vol{vol}
	\opn\link{link} \opn\star{star}\opn\lex{lex}\opn\set{set}
	\opn\gr{gr}
	
	\def\pot#1#2{#1[\kern-0.28ex[#2]\kern-0.28ex]}

	\opn\dirlim{\underrightarrow{\lim}}
	\opn\inivlim{\underleftarrow{\lim}}
	\def\Implies{\ifmmode\Longrightarrow \else
		\unskip${}\Longrightarrow{}$\ignorespaces\fi}
	\def\implies{\ifmmode\Rightarrow \else
		\unskip${}\Rightarrow{}$\ignorespaces\fi}
	\def\iff{\ifmmode\Longleftrightarrow \else
		\unskip${}\Longleftrightarrow{}$\ignorespaces\fi}

	\let\:=\colon
	\newtheorem{Theorem}{Theorem}[section]
	\newtheorem{Lemma}[Theorem]{Lemma}
	\newtheorem{Corollary}[Theorem]{Corollary}
	\newtheorem{Proposition}[Theorem]{Proposition}
	\newtheorem{Remark}[Theorem]{Remark}
	
	\newtheorem{Example}[Theorem]{Example}

	\newtheorem{Notation}[Theorem]{Notation}

	\let\epsilon\varepsilon
	\let\kappa=\varkappa	
\begin{document}
\title[] {$t$-Young Complexes and Squarefree Powers of $t$-Path Ideals}
	
\author{Francesco Navarra}
\address{Sabanci University, Faculty of Engineering and Natural Sciences, Orta Mahalle, Tuzla 34956, Istanbul, Turkey}
\email{francesco.navarra@sabanciuniv.edu}
    
\author{Ayesha Asloob Qureshi}
\address{Sabanci University, Faculty of Engineering and Natural Sciences, Orta Mahalle, Tuzla 34956, Istanbul, Turkey}
\email{aqureshi@sabanciuniv.edu, ayesha.asloob@sabanciuniv.edu}

\author {Dharm Veer}
\address{Department of Mathematics \& Statistics, Dalhousie University, 6297 Castine Way, PO
BOX 15000, Halifax, NS, Canada B3H 4R2}
\email{d.veer@dal.ca}
		
	

\thanks{The first and second authors are supported by Scientific and Technological Research Council of Turkey T\"UB\.{I}TAK under the Grant No: 124F113, and are thankful to T\"UB\.{I}TAK for their support. The first author is a member of GNSAGA-INDAM and he acknowledges its support. The last author is supported by an AARMS Postdoctoral Fellowship. 
}

\begin{abstract}
We introduce a new class of simplicial complexes, called \emph{$t$-Young complexes}, arising from Young diagrams and a positive integer~$t$. 
We show that every $t$-Young complex is either contractible or homotopy equivalent to a wedge of spheres, and a complete characterization of their vertex-decomposability is provided.
Interestingly, $t$-Young complexes naturally appear as the Alexander dual complexes of squarefree powers of $t$-path ideals of path graphs.
For this family, we use generating functions to obtain an explicit formula for homotopy type. 
As applications, we determine the projective dimension and the Krull dimension of these squarefree powers.
\end{abstract}

\maketitle

\section*{Introduction}

The study of algebraic invariants of monomial ideals through combinatorial and topological methods has been a central theme in combinatorial commutative algebra.
In this paper, we investigate the projective dimension and Krull-dimension of squarefree powers of $t$-path ideals of path graphs, combining algebraic techniques with a new topological framework based on simplicial complexes. We introduce a new class of simplicial complexes, which we call \emph{$t$-Young complexes}, and show that $t$-Young complexes exhibit nice topological behavior: each is either contractible or homotopy equivalent to a wedge of spheres. These complexes naturally arise as Alexander duals of Stanley--Reisner complexes associated with squarefree powers of $t$-path ideals. Through this connection, we obtain an explicit formula for the projective dimension of the squarefree powers of $t$-path ideals. 

The notion of $t$-path ideals of graphs, introduced in \cite{CN}, refers to monomial ideals generated by the monomials corresponding to paths of length $t-1$ in a graph $G$. In general, $t$-path ideals have been extensively studied and in the case of a rooted tree, they coincide with the facet ideal of a simplicial tree (see \cite{HT2010pathideals}). In the language of hypergraphs, simplicial trees correspond to totally balanced hypergraphs (see \cite[Theorem 3.2]{HHTZ}), that is, hypergraphs containing no special cycles. For further details, we refer the reader to \cite[Chapter 5]{B}. The combinatorial interpretation of certain algebraic invariants of $t$-path ideals of path graphs has been explored in \cite{AF, AF2, T}, and more recently, this line of research was extended to ordinary powers in \cite{BCV} and to squarefree powers in \cite{KNQ24squarefree}, where Castelnuovo–Mumford regularity was explicitly determined.

Given a monomial ideal $I$, its $k$-th \emph{squarefree power}, denoted by $I^{[k]}$, is generated by the squarefree monomials among the generators of~$I^k$.  The study of squarefree powers was formally initiated by \cite{BHZ18powersedgeideals} in the context of edge ideals of graphs, although such ideals had implicitly appeared earlier (for instance, in the proof of Proposition 2.10 in \cite{JZ}). Since then, numerous works have appeared, such as \cite{CFL, EF, EHHM, EHHM2, EH, FHH, Kamalesh1, Kamalesh2, S}. Squarefree powers provide a combinatorial shadow of the ordinary powers of~$I$, and they are often easier to analyze while still encoding important homological information. 
Another interesting aspect of the study of squarefree powers arises from their deep connection with the matching theory of simplicial complexes. Let $\Delta$ be a simplicial complex, and let $I(\Delta)$ denote its facet ideal. A \emph{matching} of $\Delta$ is a collection of pairwise disjoint facets of $\Delta$. The minimal generators of $I(\Delta)^{[k]}$ bijectively correspond to the matchings of $\Delta$ of size $k$. In particular, $I(\Delta)^{[k]} \neq 0$ if and only if $1 \leq k \leq \nu(\Delta)$, where $\nu(\Delta)$ denotes the maximum size of a matching of $\Delta$. One of our main results gives a formula for the projective dimension of the $k$-th squarefree powers of $t$-path ideals of path graphs $P_n$. 

\begin{Theorem}[Theorem \ref{Theorem: Projective dimension}]
Let $P_n$ be the path graph on $n$ vertices. Let $I_{n,t}$ and $\Gamma_{n,t}$ denote the $t$-path ideal of $P_n$ and the facet complex of $I_{n,t}$, respectively, for $1\leq t \leq n$. Then for any $1\leq k \leq \nu(\Gamma_{n,t})$, we have
        $$ 
        \mathrm{pd}\left(R/I_{n,t}^{[k]}\right)= \begin{cases}
         n-kt+1 & \text{ if }  kt \leq n \leq k(t+1); \\[0.6em]
         \displaystyle\frac{2(n-d)}{t+1}-k+1 & \text{ if } n \equiv d \mod(t+1),  \text{ with } 0\leq d \leq t-1 \text{, and } n>k(t+1);\\[0.9em]
         \displaystyle \frac{2(n+1)}{t+1}-k & \text{ if } n \equiv t \mod(t+1)  \text{ and } n> k(t+1).\\    
        \end{cases}
        $$	
\end{Theorem}

The reader may refer to Example~\ref{Exa: proj dim formula} to observe the nice pattern when $k$ and $t$ are fixed and $n$ varies.
The proof combines algebraic and topological arguments. We prove the upper bound using techniques in commutative algebra, via exact sequences. The lower bound is obtained through a topological analysis of the Alexander dual of the Stanley--Reisner complex of the ring $R/I_{n,t}^{[k]}$, which leads to the study of $t$-Young complexes. 
The Alexander dual $\Delta^\vee$ of a simplicial complex~$\Delta$ is often more convenient to understand the minimal free resolution of the Stanley--Reisner ring $K[\Delta]$ of $\Delta$, since Hochster’s formula~\cite{hochsterformula} can be expressed in terms of the link of faces in $\Delta^\vee$. There is a rich literature relating free resolutions of $K[\Delta]$ to the properties of the Alexander dual $\Delta^\vee$ (see e.g., \cite{ER98linearresolution, terai99Alexanderduality,HRW99componentwiselinear}).

The paper is organized as follows. Section~\ref{Section: Preliminaries} presents the preliminary notions. Section~\ref{Section: t-Young Complexes} is devoted to the definition and study of the $t$-Young complexes $\Delta_t^\lambda$. Let $\lambda = (\lambda_1, \lambda_2, \ldots, \lambda_r)$ with $\lambda_1\geq \lambda_2\geq \cdots\geq  \lambda_r$ be a partition represented by its Young diagram.  
We fill the diagram by placing the integer~$1$ in the bottom-left cell, and then assign entries so that they increase by~$t$ when moving upward along a column and by~$1$ when moving rightward along a row.  
The \emph{$t$-Young complex}~$\Delta_t^\lambda$ is the simplicial complex generated by all increasing sequences obtained from left to right in this filling (see Section \ref{Section: t-Young Complexes} for a formal definition).  We establish two key lemmas (Lemmas~\ref{lem:wedgeproduct1} and~\ref{lem:wedgeproduct2}) describing the homotopy type of $t$-Young complexes, and using these lemmas, we prove our first main theorem on $t$-Young complexes,

\begin{Theorem}[Theorem \ref{thm:homotopytype}]
Let $t$ be a positive integer and $\lambda = (\lambda_1, \lambda_2, \ldots, \lambda_r)$ a partition.  
Then the $t$-Young complex $\Delta_t^{\lambda}$ is either contractible or homotopy equivalent to a wedge of spheres.
\end{Theorem}

 Simplicial complexes that are contractible or homotopy equivalent to wedges of spheres play a central role in topological combinatorics; examples include shellable complexes~\cite{BW}, sequentially Cohen-Macaulay complexes~\cite{BWV09SCM} and matching and independence complexes of various graphs~\cite{Matsushita19matching, anurag20forest, kim22cyclelength3, Bayer23matching}. We refer the reader to~\cite{jonsson08book} for a comprehensive account of results establishing such homotopy types for numerous complexes arising in topological combinatorics.

We also examine when $\Delta_t^{\lambda}$ is vertex-decomposable. It turns out, when $\lambda_2 \le t$, the complex $\Delta_t^{\lambda}$ coincides with the order complex of a poset of dimension at most two, and hence is shellable by~\cite{JV25dimension2posets}, see Remark \ref{rem:2dimposet}.
In general, we obtain the following characterization of vertex-decomposability for $\Delta_t^{\lambda}$ in terms of~$\lambda$.

\begin{Theorem}[Theorem \ref{thm:t=1vd}, Theorem \ref{thm:cmdiagram}]
Let $\lambda = (\lambda_1, \lambda_2, \ldots, \lambda_r)$ be a partition.  
Then $\Delta_1^{\lambda}$ is vertex-decomposable.  
Moreover, for $t \ge 2$, the following are equivalent:
\begin{enumerate}
    \item $\Delta_t^{\lambda}$ is vertex-decomposable;
    \item $\Delta_t^{\lambda}$ is shellable;
    \item $\Delta_t^{\lambda}$ is Cohen--Macaulay;
    \item either $r=1$, or $r\ge2$ and $\lambda_2 \le t$.
\end{enumerate}
\end{Theorem}

Section~\ref{Section: Specialization to t-path} establishes a connection between combinatorial topology and commutative algebra. 
In particular, we show that the Alexander dual $\Delta_{n,t}^{[k]}$ of the Stanley--Reisner complex of $R/I_{n,t}^{[k]}$ coincides with a $t$-Young complex associated with a specific partition (see Proposition~\ref{lem:dualfacets}). 
As a consequence of Proposition~\ref{lem:dualfacets}, we observe that $t$-Young complexes arise as the Alexander duals of certain ideals generated by subsets of the minimal generators of $I_{n,t}^{[k]}$. 
This correspondence allows us to compute the Helly number of the Alexander dual of every $t$-Young complex (Corollary~\ref{cor:helly}). Finally, using a generating-function argument, we determine the homotopy type of $\Delta_{n,t}^{[k]}$ explicitly in terms of $n,k$ and $t$. 
\begin{Theorem}[\cref{thm:homotopyformula}]
    Let $k,t\geq 1$ and $n>kt$ be integers, and set $s=n-kt$. Then, the $t$-Young complex $\Delta_{n,t}^{[k]}$ is homotopy equivalent to
    \[
   \bigvee_{q=0}^{\lfloor s/(t+1)\rfloor} \bigvee_{{k+q-1 \choose q}{k \choose s-(t+1)q}}S^{s-(t-1)q-1}\,,
    \]
    where $\bigvee$ denotes the wedge sum, $S^r$ denotes the $r$-dimensional sphere, and  the empty wedge is a point.
\end{Theorem}

Section~\ref{Section: Projective dimension} is devoted to proving the formulas of the projective and Krull dimensions of $R/I_{n,t}^{[k]}$. The proof of projective dimension (Theorem \ref{Theorem: Projective dimension}) is divided into two parts, providing the lower and upper bounds, respectively in Subsections~\ref{subsection:lower bound} and~\ref{subsection:upper bound}.
As a consequence, when a $t$-Young complex arises as the Alexander dual of the Stanley--Reisner complex of the ring $R/I_{n,t}^{[k]}$, we determine its Leray number by applying our formula for the projective dimension. We conclude the paper by giving a combinatorial expression for the Krull dimension of $R/I_{n,t}^{[k]}$ (Theorem \ref{Theorem: Krull dimension}).

\subsection*{Acknowledgements} This project started when the last author visited the first and second authors at Sabancı University, Turkey, he thank them for the hospitality. The computer
algebra systems Macaulay2~\cite{M2} and SageMath~\cite{sage} provided valuable assistance in studying examples.

\section{Preliminaries}\label{Section: Preliminaries}

In this section, we recall some basic concepts concerning simplicial complexes and squarefree monomial ideals. The notation and definitions introduced here will be used throughout the paper.

\subsection*{Combinatorial Topology.}\label{SubSection: Preliminaries 1} An \textit{abstract simplicial complex} $\Delta$ on the vertex set $V(\Delta)$ is a collection of subsets of $V(\Delta)$ such that whenever $F' \in \Delta$ and $F \subseteq F'$, then $F \in \Delta$. In our definition, we allow the void complex $\emptyset$ in degenerate
cases. Notice that the void complex is different from $\{\emptyset\}$. The elements of $\Delta$ are called \textit{faces}. For any $F \in \Delta$, the \textit{dimension} of $F$, denoted by $\dim(F)$, is one less than the cardinality of $F$.  
A \textit{vertex} of $\Delta$ is a face of dimension $0$, and an \textit{edge} is a face of dimension $1$.  The dimension of $\Delta$ is given by $\max\{\dim(F) : F \in \Delta\}$.  
The maximal faces of $\Delta$ with respect to inclusion are called \textit{facets}, and the set of all facets of $\Delta$ is denoted by $\mathcal{F}(\Delta)$.

A \textit{subcollection} of $\Delta$ is a simplicial complex $\Delta'$ such that
$\mathcal{F}(\Delta')\subseteq\mathcal{F}(\Delta)$. A subcollection $\Delta'$ of $\Delta$ is said to be \textit{induced} if it contains every facet of $\Delta$ whose vertices belong to $V(\Delta')$. A \textit{subcomplex} of $\Delta$ is a simplicial complex $\Gamma$ such that $\Gamma\subseteq\Delta$. For a subset $W\subseteq V(\Delta)$, the \textit{induced subcomplex} of $\Delta$ on $W$, denoted by $\Delta_W$, is defined by $\Delta_W=\{F\in\Delta:F\subseteq W\}$. A subcomplex $\Gamma$ of $\Delta$ is said to be \textit{induced} if $\Gamma=\Delta_{V(\Gamma)}$.

Given a collection $\mathcal{F} = \{F_1, \dots, F_m\}$ of subsets of $V(\Delta)$, we denote by $\langle F_1, \dots, F_m \rangle$, or briefly $\langle \mathcal{F} \rangle$, the abstract simplicial complex consisting of all subsets of $V(\Delta)$ contained in some $F_i$, for $i = 1, \dots, m$. If $\Delta=\langle F\rangle$, that is, it is generated by only one facet, then $\Delta$ is called a \textit{simplex}. An abstract simplicial complex $\Delta$ is called \textit{pure} if all its facets have the same dimension. For a pure abstract simplicial complex, the dimension of $\Delta$ coincides with that of its facets. The \textit{Alexander dual} of $\Delta$ is defined as $\{ [n] \setminus F : F \notin \Delta \}$ and it is denoted by $\Delta^\vee$. After fixing a vertex $v \in V(\Delta)$ we can define the following abstract simplicial complexes:
\begin{itemize}
    \item the \emph{link} of $v$ in $\Delta$, as $\lk_\Delta(v) := \{\sigma \in \Delta : v \notin \sigma \text{ and } \sigma \cup \{v\} \in \Delta\}$;
    \item the \emph{deletion} of $v$ in $\Delta$, as $\del_\Delta(v) := \{\sigma \in \Delta : v \notin \sigma\}$;
    \item the \emph{star} of $v$ in $\Delta$, as $\st_\Delta(v) := \{\sigma \in \Delta : \sigma \cup \{v\} \in \Delta\}$.
\end{itemize}

We now recall the definition of vertex decomposability. This notion was first introduced by Provan and Billera in \cite{PB} for pure simplicial complexes, and later generalized by Björner and Wachs in \cite{BW97part2} to the non-pure case. Throughout this paper, we adopt \cite[Definition 11.1]{BW97part2}, although all simplicial complexes under consideration are assumed to be pure. 
A vertex $x$ of a simplicial complex $\Delta$ is called a \textit{shedding vertex} if every facet of $\del_\Delta(x)$ is also a facet of $\Delta$. A simplicial complex $\Delta$ is said to be \textit{vertex decomposable} if one of the following conditions holds:
\begin{enumerate}[(1)]
    \item $\Delta = \{\emptyset\}$, or $\Delta$ is a simplex; or
    \item there exists a shedding vertex $x \in V(\Delta)$ such that both $\lk_\Delta(x)$ and $\del_\Delta(x)$ are vertex decomposable.
\end{enumerate}

We now provide the definition of a large class of vertex-decomposable simplicial complexes. A pure simplicial complex $\Delta$ is called \textit{shellable} if the facets of $\Delta$ can be ordered as $F_1,\dots,F_m$ in such a way that $\langle F_1,\dots,F_{i-1}\rangle \cap \langle F_i\rangle$ is generated by a non-empty set of maximal proper faces of $F_i$, for all $i\in \{2,\dots,m\}$. It is known from \cite[Theorem 11.3]{BW97part2} that every vertex-decomposable simplicial complex is shellable.

The \emph{join} of two abstract simplicial complexes $\Delta_1$ and $\Delta_2$, denoted by $\Delta_1 * \Delta_2$, is the abstract simplicial complex whose simplices are disjoint unions of simplices of $\Delta_1$ and $\Delta_2$.  
For $v \notin V(\Delta)$, the \emph{cone} on $\Delta$ with apex $v$, denoted by $C_v(\Delta)$, is defined as
\[
C_v(\Delta) = \Delta * \langle \{v\} \rangle.
\]
For distinct vertices $v, w \notin V(\Delta)$, the \emph{suspension} of $\Delta$, denoted by $\Sigma(\Delta)$, is defined as
\[
\Sigma(\Delta) = \Delta * \langle \{v\}, \{w\} \rangle.
\]
We denote by $S^n$ the $n$-dimensional sphere, for $n\geq 0$. Moreover, we adopt the convention that $S^{-1}=\{\emptyset\}$; in particular, $\Sigma(\{\emptyset\}) =S^0$.

In the literature, it is customary to denote by $\lvert \Delta \rvert$ the \textit{geometric realization} of an abstract simplicial complex $\Delta$. With a slight abuse of notation, and in order to avoid cumbersome distinctions, we use the same symbol $\Delta$ to denote both an abstract simplicial complex and its geometric realization, referring to either simply as a \textit{simplicial complex}; the intended meaning will always be clear from the context. 

The following lemma, which will play a crucial role in our discussion, follows from \cite[Example~0.14]{hatcher}. We recall that a simplicial complex $\Delta$ is said to be \textit{contractible} if it is homotopy equivalent to a single point; moreover, if $\Delta'$ is a subcomplex of $\Delta$, we say that $\Delta'$ is \emph{contractible in} $\Delta$ if the inclusion map $\Delta'\hookrightarrow\Delta$ is homotopic to a constant map. For a more detailed account of the algebraic topology underlying these notions, we refer the reader to \cite{kozlov2007combinatorial}.

\begin{Lemma}\label{lem:cone}
    Let $\Delta$ be a simplicial complex and $v \in V(\Delta)$.  
    If $\lk_\Delta(v)$ is contractible in $\del_\Delta(v)$, then
    \[
    \Delta \simeq \del_\Delta(v) \vee \Sigma(\lk_\Delta(v)).
    \]
    Here, $\simeq$ denotes the homotopy equivalence, and $\vee$ the wedge sum. 
\end{Lemma}

After introducing the necessary combinatorial and topological notions, we now turn to the commutative algebra associated with simplicial complexes. 

\subsection*{Commutative Algebra.}\label{SubSection: Preliminaries 2} Let $\Delta$ be a simplicial complex on $[n] := \{1, \dots, n\}$ and $R = K[x_1, \dots, x_n]$ be a polynomial ring in $n$ variables over a field $K$. A \textit{non-face} of $\Delta$ is a subset $F \subseteq V(\Delta)$ such that $F \notin \Delta$, and we denote by $\mathcal{N}(\Delta)$ the set of minimal non-faces of $\Delta$. For any subset $F = \{i_1, \dots, i_r\} \subseteq [n]$, we associate the squarefree monomial $x_F = x_{i_1} \cdots x_{i_r} \in R$. The \textit{facet ideal} of $\Delta$ is the ideal $I(\Delta)$ of $R$ generated by the squarefree monomials $x_F$ such that $F$ is a facet of $\Delta$. The \textit{Stanley-Reisner ideal} of $\Delta$ is the ideal $I_\Delta$ of $R$ generated by those squarefree monomials $x_F$ with $F \notin \Delta$, that is, $I_\Delta = (x_F : F \in \mathcal{N}(\Delta))$. Moreover, the quotient ring $R/I_{\Delta}$ is called the \textit{Stanley-Reisner ring} of $\Delta$ and is denoted by $K[\Delta]$. We recall that $\Delta$ is said to be \textit{Cohen-Macaulay} over $K$ if $\widetilde{H}_i(\lk_\Delta(\sigma), K)= 0$ for all $i< \dim(\lk_\Delta(\sigma))$ for every face $\sigma$ of $\Delta$, where $\lk_\Delta(\sigma)=\{ \tau \in \Delta: \tau \cap \sigma = \emptyset, \tau \cup \sigma \in \Delta\}$. Reisner~\cite[Theorem\ 1]{R} proved that $\Delta$ is Cohen-Macaulay over $K$ if and only if the Stanley-Reisner ring associated to $\Delta$ (over $K$) is Cohen-Macaulay.
From \cite[Theorem 5.1.13]{BH}, it is known that every shellable simplicial complex is Cohen-Macaulay.

From an algebraic combinatorial point of view, simplicial complexes play an essential role in the study of squarefree monomial ideals, that is, ideals generated by monomials that are products of distinct variables. In fact, every squarefree monomial ideal can be regarded as the Stanley–Reisner ideal of a simplicial complex. In this paper, we will focus on the so-called squarefree powers, a notion of power introduced in \cite{BHZ18powersedgeideals}. Given a squarefree monomial ideal $I$ in $R=K[x_1, \ldots, x_n]$, the $k$-th squarefree power of $I$ is defined to be the ideal generated by the squarefree elements of $G(I^k)$, and it is denoted by $I^{[k]}$. A {\em matching} of $\Delta$ is a set of pairwise disjoint facets of $\Delta$. A matching consisting of $k$ facets is referred to as a \textit{$k$-matching}. A $k$-matching is called {\em maximal} if $\Delta$ does not admit any $(k+1)$-matching. The {\em matching number} of $\Delta$ is the size of a maximum matching of $\Delta$ and is denoted by $\nu(\Delta)$. When $I$ is regarded as the facet ideal of a simplicial complex $\Delta$, the monomial generators of $I^{[k]}$ correspond to the $k$-matchings of $\Delta$, that is,
\[
I^{[k]} = (x_{F_1}\cdots x_{F_k} : \{F_1, \ldots, F_k\} \text{ is a $k$-matching of } \Delta).
\]
It is straightforward to check that $I^{[k]} \neq 0$ if and only if $1 \leq k \leq \nu(\Delta)$. In particular, $I^{[\nu(\Delta)]}$ is the highest non-vanishing squarefree power of $I$.\\
A further kind of matching can be defined for a simplicial complex, namely the \textit{restricted matching}. Let $F, G \in \mathcal{F}(\Delta)$. Then $F$ and $G$ form a {\em gap} in $\Delta$ if $F \cap G = \emptyset$ and the induced subcollection on the vertex set $F \cup G$ is $\langle F, G \rangle$. A matching $M$ of $\Delta$ is called a {\em restricted matching} if there exists a facet in $M$ forming a gap with every other facet in $M$. The maximal size of a restricted matching of $\Delta$ is denoted by $\nu_0(\Delta)$ and is called the \textit{restricted matching number} of $\Delta$.

Recall that a $1$-dimensional simplicial complex is called a \textit{graph}. A \textit{path} of length $r-1$ in a graph $G$ is a sequence of $r$ distinct vertices $v_1, \dots, v_r$ such that $\{v_i, v_{i+1}\}$ is an edge of $G$ for all $i = 1, \dots, r-1$. In \cite{CN}, the authors introduced the $t$-path ideals of graphs, for $t$ positive integer. Specifically, if $G$ is a graph, the \emph{$t$-path ideal} of $G$ is defined by
\[
I_t(G) = (x_{i_1} \cdots x_{i_t} : \{i_1, \ldots, i_t\} \text{ forms a path on $t$ vertices in } G).
\]
Observe that $I_1(G)=(x_1,\dots,x_n)$, if $G$ is a graph on $[n]$, and $I_2(G)$ coincides with the so-called \textit{edge ideal} of $G$ (see \cite[Subsection 9.1]{HHmonomialideals}). A graph on the vertex set $[n]$ forming a path of length $n - 1$ is called a \textit{path graph}, and is denoted by $P_n$. We conclude the preliminaries by introducing the notation that will be used throughout the paper, starting from Section \ref{Section: Specialization to t-path}.

\begin{Notation}\label{Notation}
Let $P_n$ be the path graph on $n$ vertices. We denote by $I_{n,t}$ the $t$-path ideal of $P_n$, for $1\leq t \leq n$, and by $I_{n,t}^{[k]}$ its $k$-th squarefree power, for $k \geq 1$. Then we set:
\begin{itemize}
    \item $\Gamma_{n,t}$ to denote the facet complex of $I_{n,t}$, that is, $I(\Gamma_{n,t})=I_{n,t}.$
    \item $\Sigma_{n,t}^{[k]}$ to denote the Stanley–Reisner complex of $I_{n,t}^{[k]}$, that is, $I_{\Sigma_{n,t}^{[k]}}=I_{n,t}^{[k]}.$ 
    \item $\Delta_{n,t}^{[k]}$ to denote the Alexander dual of $\Sigma_{n,t}^{[k]}$, that is, $\Delta_{n,t}^{[k]} = (\Sigma_{n,t}^{[k]})^\vee$.
\end{itemize}
\end{Notation}
 
In addition to the above notation, we set $\Delta_{n,t}^{[0]}$ as the simplex $\langle\{1,2,\ldots,n\}\rangle$.

\section{$t$-Young complexes} \label{Section: t-Young Complexes}

In this section, we introduce the main combinatorial object studied throughout the paper, which we call the \textit{$t$-Young complex} associated with a partition $\lambda$. This simplicial complex is obtained by filling the Young diagram of $\lambda$ according to a fixed parameter $t$ so that the entries increase by $t$ along the first column, and each row is subsequently filled from left to right with consecutive integers increasing by one. 

Let $\lambda = (\lambda_1, \lambda_2,  \ldots,  \lambda_r)$  with $\lambda_1 \geq \lambda_2 \geq  \ldots \geq  \lambda_r \geq 1$ be a partition,  identified with its Young diagram: a top justified shape consisting of $r$ rows of boxes of lengths $\lambda_1, \ldots,\lambda_r$ with the row sizes weakly decreasing. For a positive integer $t$, and for $1\leq j\leq r$, let 
\[
\MS_j  = \{(r-j)t+ i : 1\leq i \leq \lambda_j\}.
\]

Fill the boxes of the $j$-th row in the Young diagram $\lambda$ by entries of $\MS_j$ from left to right in increasing order. Then, we define the following $(\lambda_1-1)$-dimensional pure simplicial complex generated by
\[
 \Big\langle
\{\, a_{1}, a_{2}, \ldots, a_{\lambda_{1}} \,\}
\ \Big|\ 
a_{1} < a_{2} < \cdots < a_{\lambda_{1}},
\ 
a_{j} \text{ lies in the $j$-th column of } \lambda 
\text{ for all } 1 \le j \le \lambda_{1}
\Big\rangle,
\]
which we call \textit{$t$-Young complex} and denote by $\Delta_{t}^{\lambda}$. 

Although partitions are assumed to have positive parts, the zero tuple may occur in some of the recursive formulas below. Therefore, for $\mu=(0,\ldots,0)$, we adopt the convention $\Delta_t^\mu=\{\emptyset\}.$

\begin{Example} \label{eg:young}
\rm
    The Young diagram corresponding to the partition $\lambda = (5,4,2)$ and $t = 3$ is shown in Table \ref{tab:exa}. In this case, the facets of $\Delta_{t}^\lambda$ are the following: $\{1,2,6,7,11\}$, $\{1,2,6,10,11\}$,  $\{1,2,9,10,11\}$, $\{1,5,6,7,11\}$, $\{1,5,6,10,11\}$, $\{1,5,9,10,11\}$, $\{1,8,9,10,11\}$, $\{4,5,6,7,11\}$, $\{4,5,6,10,11\}$, $\{4,5,9,10,11\}$, $\{4,8,9,10,11\}$ and $\{7,8,9,10,11\}$.

\begin{table}[ht]
\centering
\begin{ytableau}
7 & 8 & 9 & 10 & 11\\
4 & 5 & 6 & 7 & \none \\ 
1 & 2 & \none & \none & \none \\ 
\end{ytableau}
\caption{Young diagram filled by $\MS_1$, $\MS_2$ and $\MS_3$.}
\label{tab:exa}
\end{table}
\end{Example}

In the following remark, we see that $\Delta_{t}^{\lambda}$ is the order complex of a poset for $\lambda_2\leq t$. We refer the reader to \cite{JV25dimension2posets} for the terminology and definitions used in the remark.

\begin{Remark}\label{rem:2dimposet}
Let $t$ be a positive integer and $\lambda = (\lambda_1, \lambda_2,  \ldots,  \lambda_r)$ be a partition such that $\lambda_2\leq t$. Define a poset $P_t^{\lambda}$ on the set $\cup_{i\in [r]} \MS_i$ with the following cover relations: for $x,y \in P_t^{\lambda},$ we have $x<y$ in $P_t^{\lambda}$ if and only if $x<y$ in $\NN$, and $x$ and $y$ are in $i$-th and $(i+1)$-th column, respectively for some $1\leq i <\lambda_1$. It is not hard to observe that $\Delta_t^{\lambda}$ is the order complex of the poset $P_t^{\lambda}$.

Let $\sigma$ be a permutation with the following entries: entries of first column of $\lambda$ in descending order followed by entries of second column of $\lambda$ in descending order and so on. Then $P_t^{\lambda}$ is the intersection of the identity permutation of $\cup_{i\in [r]} \MS_i$ and $\sigma$.
Therefore, $P_t^{\lambda}$ is a poset of dimension $2$ if $r\geq 2$ and of dimension $1$ if $r=1$. It follows from \cite[Theorem\ 1]{JV25dimension2posets} that $\Delta_t^{\lambda}$ is shellable.

\end{Remark}

In the sequel, we show that $t$-Young complexes are either contractible or homotopy equivalent to a wedge of spheres. To this end, we establish two lemmas describing the homotopy type of $t$-Young complexes when the first two rows of the Young diagram have the same length. We begin with the case where this common length is at most $t$.

\begin{Lemma}\label{lem:wedgeproduct1}
Let $t$ be a positive integer and $\lambda = (\lambda_1, \lambda_2,  \ldots,  \lambda_r)$  with $\lambda_1 =\lambda_2$ be a partition.
 If $\lambda_1\leq t,$ then 
    \[
    \Delta_t^{\lambda}\simeq \Delta_t^{\lambda'} \vee \Sigma(\Delta_t^{\lambda''}),
    \] 
where $\lambda' = (\lambda_2, \lambda_3,  \ldots,  \lambda_r),$ and $\lambda'' = (\min\{\lambda_1-1, \lambda_2\},\min\{\lambda_1-1, \lambda_3\}, \ldots,\min\{\lambda_1-1, \lambda_r\})$.
\end{Lemma}

\begin{proof}
We first consider the trivial case $\lambda_1=1$. Since $\lambda_1=\lambda_2$ and $\lambda$ is a partition, we have $\lambda=(1,\ldots,1)$. Moreover, $\lambda''=(0,\ldots,0)$, and therefore
$\Delta^{\lambda''}_t=\{\emptyset\}$ and $\Sigma(\Delta^{\lambda''}_t)=S^0$. Consequently, $\Delta^\lambda_t \simeq \Delta^{\lambda'}_t \vee \Sigma(\Delta^{\lambda''}_t)$.

Thus, we may assume that $\lambda_1\geq2$. Since $\lambda_1=\lambda_2$, the top two rows of the Young diagram are as shown in Table \ref{tab: top two rows 1}.

  \begin{table}[h]
\begin{tabular}{ |c|c|c|c | } 
 \hline
 $(r-1) t+ 1$ & $(r-1) t+ 2$ & $\cdots$ & $(r-1) t+ \lambda_1$\\
 \hline
 $(r-2) t+1$ & $(r-2) t+ 2$ & $\cdots $ & $(r-2) t+ \lambda_2$ \\ 
 \hline
\end{tabular}
\caption{Top two rows of Young diagram, when $\lambda_1 =\lambda_2$ and $\lambda_1\leq t$.}
\label{tab: top two rows 1}
\end{table}

Set $v_i = (r-1) t+ i$ for each $i=1, \ldots, \lambda_1$. Because $\lambda_1\leq t,$ these vertices appear only in the top row of the Young diagram. 
We proceed by deleting the vertices $v_1, \ldots, v_{\lambda_1}$ one by one and analyzing the corresponding links. Set $\Delta_0=\Delta_t^\lambda$ and define recursively $\Delta_i=\del_{\Delta_{i-1}}(v_i)$ for $1\leq i\leq\lambda_1$.

\medskip
\noindent\textbf{Step 1: Deletion of $v_1$.} 
By the definition of the facets of $\Delta_t^{\lambda}$, it follows that the only facet containing $v_1$ is $\{v_1, v_2, \ldots, v_{\lambda_1}\}$. Thus, $\lk_{\Delta_t^{\lambda}}(v_1) = \langle\{v_2,\ldots,v_{\lambda_1}\}\rangle,$ which is a cone. Therefore, it is contractible in $\Delta_1$. Hence,  
\begin{equation}\label{eq:v1case1}
    \Delta_t^{\lambda} \simeq \Delta_1 \vee \Sigma(\lk_{\Delta_t^{\lambda}}(v_1)) \simeq \Delta_1,
\end{equation}
where the first homotopy equivalence follows from Lemma~\ref{lem:cone}, and the second follows from the contractibility of $\lk_{\Delta_t^{\lambda}}(v_1)$.

\medskip
\noindent\textbf{Step 2: Deletion of $v_2, \ldots, v_{\lambda_1 - 1}$.} If $\lambda_1=2$, there is nothing to prove in this step. Assume that $\lambda_1\geq3$, fix $i \in \{1, \ldots, \lambda_1-2\}$, and consider the complex $\Delta_i$. We now examine $\lk_{\Delta_i}(v_{i+1})$, with reference to Lemma~\ref{lem:cone}, and show that it is a cone. 

\begin{table}[h]
\centering
\begin{tabular}{ |c|c|c|c | c |c | c | }
\hline
& $\cdots$ & & $v_{i+1}=(r-1) t+ i+1$ &\cellcolor{lightgray} $v_{i+2}$ &\cellcolor{lightgray} $\cdots$ &\cellcolor{lightgray} $v_{\lambda_1}=(r-1) t+ \lambda_1$ \\
\hline
\cellcolor{lightgray} $v_1-t$ & \cellcolor{lightgray} $\cdots$ & \cellcolor{lightgray} $v_{i}-t$ & $v_{i+1}-t$ & $v_{i+2}-t$ & $\cdots$ & $v_{\lambda_1}-t$ \\
\hline
\cellcolor{lightgray} $\vdots$ & \cellcolor{lightgray} $\cdots$ & \cellcolor{lightgray} $\vdots$ & $\vdots$ & $\vdots$ & $\cdots$ & $\vdots$ \\
\end{tabular}
\caption{In $\lk_{\Delta_i}(v_{i+1})$, the faces are defined using entries from the shaded region.}
\label{fig:link}
\end{table}

Let $F$ be a facet of $\lk_{\Delta_i}(v_{i+1})$. Then $v_{i+1} \notin F$ and $F \cup \{v_{i+1}\} \in \Delta_i$. Since $\Delta_i\subset \Delta_t^{\lambda}$ and by the definition of the faces of $\Delta_t^{\lambda}$, we have $F \cup \{v_{i+1}\} = \{a_1, \ldots, a_{i}, v_{i+1}, a_{i+2}, \ldots, a_{\lambda_1}\}$,
where $a_1 < \cdots < a_{i} < v_{i+1} < a_{i+2} < \cdots < a_{\lambda_1}$ and each $a_j$ lies in the $j$-th column of $\lambda$ for all $1 \le j \le \lambda_1$ with $j\neq i+1$. Since $\Delta_i = \del_{\Delta_{i-1}}(v_i)$, no face of $\Delta_i$ contains $v_1, \ldots, v_i$. Hence, for all $1 \le j \le i$, the vertex $a_j$ lies in the shaded grey region on the left in Table \ref{fig:link}. Moreover, observe that we necessarily have $a_{j} = v_{j}$, for all $i+2 \le j \le \lambda_1 $. Therefore, $\lk_{\Delta_i}(v_{i+1}) = \Delta_t^{\mu} * \langle {v_{i+2}, \ldots, v_{\lambda_1}} \rangle$, where $\mu = (\min\{\lambda_2, i\}, \ldots, \min\{\lambda_r, i\})$. In particular, $\lk_{\Delta_i}(v_{i+1})$ is a cone, as desired.

Consequently, by the same argument as before, 
\begin{equation}\label{eq:vicase1}
    \Delta_i \simeq \Delta_{i+1}, \qforall  1 \leq i \leq {\lambda_1}-2.
\end{equation}

\medskip
\noindent\textbf{Step 3: Deletion of $v_{\lambda_1}$.} 
Finally, we need to resolve $\Delta_{{\lambda_1}-1}$ (see Table \ref{fig:deletion}). 

\begin{table}[h]
    \centering
\begin{tabular}{ |c|c|c|c | } 
 \hline
  & $\cdots$ &   & $v_{\lambda_1}=(r-1) t+ \lambda_1$\\
 \hline
\cellcolor{lightgray} $v_1-t$ & \cellcolor{lightgray} $\cdots$ & \cellcolor{lightgray} $v_{\lambda_1-1}-t$  &  $v_{\lambda_1}-t$ \\ 
 \hline
 \cellcolor{lightgray} $\vdots$  & \cellcolor{lightgray} $\cdots$ & \cellcolor{lightgray} $\vdots$ & $\vdots$ \\ 
\end{tabular}
    \caption{In the $\lk_{\Delta_{\lambda_1-1}}(v_{\lambda_1})$, the faces are defined by considering only the entries from the light gray (shaded) region. Meanwhile, the complex $\Delta_{\lambda_1}$ is formed using all entries below the top row.}
    \label{fig:deletion}
\end{table}

In $\Delta_{{\lambda_1}-1}$, all entries in the top row have been removed except for the final entry, $v_{\lambda_1}$. After deleting $v_{\lambda_1}$, the resulting complex $\Delta_{\lambda_1}$ is precisely $\Delta_t^{\lambda'}$, where $\lambda' = (\lambda_2, \lambda_3,  \ldots,  \lambda_r)$, since the entire top row has now been removed.\\
Next, we consider the link of $v_{\lambda_1}$ in $\Delta_{\lambda_1-1}$. By definition, the facets of $\lk_{\Delta_{\lambda_1-1}}(v_{\lambda_1})$ are of the form $F\setminus \{v_{\lambda_1}\}$ where $F$ is a facet of $\Delta_{\lambda_1-1}$ containing $v_{\lambda_1}$. 

Hence, $\lk_{\Delta_{\lambda_1-1}}(v_{\lambda_1}) = \Delta_t^{\lambda''}$, where $\lambda'' = (\mu_2,\mu_3, \ldots,\mu_r)$ with $\mu_i = \min\{\lambda_1-1, \lambda_i\}$ for all $2\leq i \leq r.$

Since $\lambda_2= \lambda_1\leq t,$ the vertex $v_{\lambda_1}-t$ in $\Delta_{\lambda_1}$ (that is, in $ \Delta_t^{\lambda'})$ appears only in the top cell of the last column. Moreover $v_{\lambda_1}-t \geq v$ for every $v \in \lk_{\Delta_{\lambda_1-1}}(v_{\lambda_1})$.  
Hence, for every facet $G$ of $\lk_{\Delta_{\lambda_1-1}}(v_{\lambda_1})$, the set $G \cup \{v_{\lambda_1}-t\}$ is a facet of $\Delta_{\lambda_1}$.  
It follows that $\lk_{\Delta_{\lambda_1-1}}(v_{\lambda_1}) \subseteq \st_{\Delta_{\lambda_1}}(v_{\lambda_1}-t)$, and hence $\lk_{\Delta_{\lambda_1-1}}(v_{\lambda_1})$ is contractible in $\Delta_{\lambda_1}$. Applying Lemma \ref{lem:cone} once again yields
\begin{equation}\label{eq:vrcase1}
\Delta_{\lambda_1-1} \simeq \Delta_{\lambda_1} \vee \Sigma(\lk_{\Delta_{\lambda_1-1}}(v_{\lambda_1})) =\Delta_t^{\lambda'} \vee \Sigma(\Delta_t^{\lambda''}). 
\end{equation}
Putting Equations (\ref{eq:v1case1}), (\ref{eq:vicase1}) and (\ref{eq:vrcase1}) together, we obtain the desired formula.
\end{proof}

Next we consider the case when $\lambda_1> t$.


\begin{Lemma}\label{lem:wedgeproduct2}
Let $t$ be a positive integer and $\lambda = (\lambda_1, \lambda_2,  \ldots,  \lambda_r)$  with $\lambda_1 =\lambda_2$ be a partition.
 If $\lambda_1> t,$ then 
    \[
    \Delta_t^{\lambda}\simeq \Delta_t^{\lambda'} \vee \Sigma(\Delta_t^{\lambda''}) \vee \Sigma^2(\Delta_t^{\lambda'''}),
    \] 
where 
\begin{itemize}
    \item $\lambda' = (\lambda_2, \lambda_3,  \ldots,  \lambda_r)$,
    \item $\lambda'' = (\min\{\lambda_1-1, \lambda_2\},\min\{\lambda_1-1, \lambda_3\}, \ldots,\min\{\lambda_1-1, \lambda_r\})$,
    \item $\lambda''' = (\min\{\lambda_1-t-1, \lambda_1\},\min\{\lambda_1-t-1, \lambda_2\}, \ldots,\min\{\lambda_1-t-1, \lambda_r\})$
\end{itemize} 
\end{Lemma}

\begin{proof}
Let $n = (r-1) t+ \lambda_1$. Since $\lambda_1>t$ and $\lambda_1=\lambda_2$, the top two rows of the Young diagram $\lambda$ are as in Table \ref{tab: top two rows 2}.

\begin{table}[h]
    \centering
   \begin{tabular}{ |c|c|c|c |c |c |c |} 
 \hline
 $(r-1) t+ 1$  & $\cdots$ & $n-t$ & $v_1=n-t+1$ & $\cdots$ & $v_{t-1}=n-1$ & $v_t=n$\\
 \hline
 $(r-2) t+1$ & $\cdots$ & $n-2t$ &  $n-2t+1$ & $\cdots $  & $n-t-1$ & $n-t$ \\ 
 \hline
\end{tabular}
    \caption{The two top rows of the Young diagram, when $\lambda_1=\lambda_2$ and $\lambda_1>t$.}
    \label{tab: top two rows 2}
\end{table}

For each $i\in \{1, \ldots, t\}$, set $v_i = n- t+ i$, and observe that each vertex $v_i$ only appears in the $(\lambda_1-t+i)$-th column of the top row.
Define $\Delta_0=\Delta_t^\lambda$ and $\Delta_i=\del_{\Delta_{i-1}}(v_i)$ for $1\leq i\leq t$. We proceed again by deleting these vertices one by one and analyzing their links (For $t=1$ one must skip Steps 1–2 and go straight to Step 3 with $\Delta_{0}=\Delta$).
The argument is similar to Lemma~\ref{lem:wedgeproduct1}, except that we need to work a bit more to resolve $\lk_{\Delta_{t-1}}(v_t)$. 

\medskip
\noindent\textbf{Step 1: Deletion of $v_1$.}
We first examine $\lk_{\Delta_t^{\lambda}}(v_1)$. Unlike in Lemma~\ref{lem:wedgeproduct1}, where the vertex $v_1$ occupied the first column of the top row, here $v_1$ lies in the $(\lambda_1 - t + 1)$-st column of the same row. As in Lemma~\ref{lem:wedgeproduct1}, we have $\lk_{\Delta_t^{\lambda}}(v_1) = \Delta_t^{\mu} * \langle\{v_2,\ldots,v_t\}\rangle$, where $\mu = (\min\{\lambda_1-t, \lambda_1\},\min\{\lambda_1-t, \lambda_2\}, \ldots,\min\{\lambda_1-t, \lambda_r\})$ (see Table \ref{fig:linkv_1}).

\begin{table}[h]
    \centering
\begin{tabular}{ |c|c|c|c |c |c |c |} 
 \hline
 \cellcolor{lightgray} $(r-1) t+ 1$  & \cellcolor{lightgray} $\cdots$ & \cellcolor{lightgray} $n-t$ & $v_1=n-t+1$ & \cellcolor{lightgray}$\cdots$ & \cellcolor{lightgray} $v_{t-1}=n-1$ &\cellcolor{lightgray} $v_t=n$\\
 \hline
 \cellcolor{lightgray} $(r-2) t+1$ & \cellcolor{lightgray} $\cdots$ & \cellcolor{lightgray} $n-2t$ &  $n-2t+1$ & $\cdots $  & $n-t-1$ & $n-t$ \\ 
 \hline
  \cellcolor{lightgray} $\vdots$  & \cellcolor{lightgray} $\cdots$ &  \cellcolor{lightgray} $\vdots$  &  $\vdots$  & $\cdots$ & $\vdots$ & $\vdots$ \\ 
\end{tabular}.
    \caption{In $\lk_{\Delta_t^{\lambda}}(v_1)$, the faces are defined by considering only the entries from the shaded region.}
    \label{fig:linkv_1}
\end{table}

Since $\lk_{\Delta_t^{\lambda}}(v_1)\subset \Delta_1$ is a cone, it is contractible in $\Delta_1$. Therefore, 
\[
\Delta_t^{\lambda} \simeq \Delta_1 \vee \Sigma(\lk_{\Delta_t^{\lambda}}(v_1)) \simeq \Delta_1,
\]
where the first homotopy equivalence follows from Lemma \ref{lem:cone}, and the second follows from the contractibility of $\lk_{\Delta_t^{\lambda}}(v_1)$.

\medskip
\noindent\textbf{Step 2: Deletion of $v_2, \ldots, v_{t - 1}$.}
Fix $i \in \{1, \ldots, t-2\}$, and consider the complex $\Delta_i$. Since $v_1, \ldots, v_i$ have been deleted, no face in $\Delta_i$ contain these vertices. Using an argument similar to that for $v_1$ in Step~1 and for $v_{i+1}$ in Lemma~\ref{lem:wedgeproduct1}, we get that  $\lk_{\Delta_i}(v_{i+1}) =   K_{i+1} * \langle\{v_{i+2}, \ldots, v_t\}\rangle$, where $K_{i+1} = \Delta_t^{\mu} \cup \Delta_t^{\mu_i}$ is a subcomplex of $\Delta_{i+1}$ for $\mu$ as in Step 1 and $\mu_i = (\min\{\lambda_1-t+i, \lambda_2\},\min\{\lambda_1-t+i, \lambda_3\}, \ldots,\min\{\lambda_1-t+i, \lambda_r\})$ (see Table \ref{fig:linkv_i}).

\begin{table}[h]
    \centering
\begin{tabular}{ |c|c|c|c |c |c |c |c |} 
 \hline
 \cellcolor{lightgray} $\cdots$ & \cellcolor{lightgray} $n-t$ &  & $\cdots$ &   & $v_{i+1} = n-t+i+1$ & \cellcolor{lightgray} $\cdots$ &  \cellcolor{lightgray} $v_t = n$\\
 \hline
 \cellcolor{lightgray} $\cdots$ & \cellcolor{lightgray} $n-2t$ & \cellcolor{lightgray} $n-2t+1$ & \cellcolor{lightgray} $\cdots $  &  \cellcolor{lightgray} $n-2t+i$ & $n-2t+i+1$ & $\cdots$ &  $n-t$ \\ 
 \hline
 \cellcolor{lightgray} $\cdots$ &  \cellcolor{lightgray} $\vdots$  & \cellcolor{lightgray} $\vdots$  & \cellcolor{lightgray} $\cdots$ &   \cellcolor{lightgray} $\vdots$  & $\vdots$  & $\cdots$ &  $\vdots$ \\ 
\end{tabular}.
    \caption{In $\lk_{\Delta_{i}}(v_{i+1})$, the faces are defined by considering only the entries from the shaded region.}
    \label{fig:linkv_i}
\end{table}

Proceeding as in Step~2 of Lemma~\ref{lem:wedgeproduct1}, we conclude that 
\[
\Delta_i \simeq \Delta_{i+1} \qforall 1\leq i \leq t-2.
\] 

\medskip
\noindent\textbf{Step 3: Deletion of $v_t$ in $\Delta_{t - 1}$.}
First, we claim that $\Delta_t = \Delta_t^{\lambda'},$ where $\lambda' = (\lambda_2, \lambda_3,  \ldots,  \lambda_r)$. Clearly, $\Delta_t^{\lambda'} \subseteq \Delta_t$ (see Table \ref{fig:Delta_t}).

\begin{table}[h]
    \centering
\begin{tabular}{ |c|c|c|c |c |c |c |} 
 \hline
 \cellcolor{lightgray} $(r-1) t+ 1$  & \cellcolor{lightgray} $\cdots$ & \cellcolor{lightgray} $n-t$ &  & $\cdots$ &  & \\
 \hline
 \cellcolor{lightgray} $(r-2) t+1$ & \cellcolor{lightgray} $\cdots$ & \cellcolor{lightgray} $n-2t$ & \cellcolor{lightgray}  $n-2t+1$ &\cellcolor{lightgray} $\cdots $  &\cellcolor{lightgray} $n-t-1$ &\cellcolor{lightgray} $n-t$ \\ 
 \hline
  \cellcolor{lightgray} $\vdots$  & \cellcolor{lightgray} $\cdots$ &  \cellcolor{lightgray} $\vdots$  & \cellcolor{lightgray} $\vdots$  & \cellcolor{lightgray} $\cdots$ & \cellcolor{lightgray} $\vdots$ &\cellcolor{lightgray} $\vdots$ \\ 
\end{tabular}.
    \caption{$\Delta_t$}
    \label{fig:Delta_t}
\end{table}

After deleting the vertices $v_1,\ldots, v_t$, the remaining $\lambda_1-t$ vertices in the top row are $(r -1)t+1, \ldots, n-t$ which coincide with the last $\lambda_1-t$ vertices of the second row. 

Now let $F$ be any face of $\Delta_t$ containing one of the surviving top-row vertices, say $(r-1)t+j \in \MS_1$ for some $1\leq j\leq \lambda_1 -t$. The entries of $F$ that are greater than $(r-1)t+j$ correspond to vertices lying to the right of it in the top row. We may reinterpret the vertex $(r-1)t+j$ as the entry in the $(j+t)$-th cell of the second row from the top. Similarly, if $F$ contains the vertex in the $(j+t)$-th cell of the top column, we may view it as coming from the $(j+2t)$-th cell of the second row, and so on.  
Since there are no top-row vertices in the last $t$ columns, $F$ can be realized as a face in the simplicial complex obtained by deleting the top row. Hence, $F\in \Delta_t^{\lambda'},$ which proves the claim.

Next, we analyze $\lk_{\Delta_{t-1}}(v_t)$. Since $v_t$ belongs to the rightmost column (the $\lambda_1$-th column) prior to deletion, the complex $\lk_{\Delta_{t-1}}(v_t)$ is supported on the first $\lambda_1-1$ columns of $\Delta_{t-1}$.
The portion of the diagram that generates $\lk_{\Delta_{t-1}}(v_t)$ is displayed in Table \ref{tab:link vt}.

\begin{table}[h]
\begin{tabular}{ |c|c|c|c |c |c |c|} 
 \hline
 $(r-1) t+ 1$  & $\cdots$ & $n-t-1$  & $n-t$ &  & &  \\
 \hline
 $(r-2) t+1$ & $\cdots$ &  $n-2t-1$ & $n-2t$ &  $n-2t+1$ & $\cdots $  & $n-t-1$ \\ 
 \hline
  $\vdots$  & $\vdots$& $\vdots$ &  $\vdots$ &  $\vdots$ & $\vdots$  & $\vdots$\\ 
\end{tabular}
  \caption{$\lk_{\Delta_{t-1}}(v_t)$}
    \label{tab:link vt}
\end{table} 

Note that $\mathrm{lk}_{\Delta_{t-1}}(v_t)=\Delta^{\lambda''}_t\cup\Delta^{\mu}_t$, with $\mu$ as in Step 1: a face containing a top-row vertex cannot extend beyond column $\lambda_1-t$, since $v_1,\dots,v_{t-1}$ have been deleted. In particular, $\mathrm{lk}_{\Delta_{t-1}}(v_t)$ need not be pure. Let $G$ be a facet of $\mathrm{lk}_{\Delta_{t-1}}(v_t)$. If $G\in\Delta^{\lambda''}_t$, then $G\cup\{n-t\}$ is a facet of $\Delta^{\lambda'}_t$, because $n-t$ appears only in the last cell of the top row of $\lambda'$ and $n-t>v$ for all $v\in G$. Otherwise $G$ is a facet of $\Delta^{\mu}_t$; then $n-t\in G$ and $G\in\Delta_t=\Delta^{\lambda'}_t$. In both cases $G\in\mathrm{st}_{\Delta_t}(n-t)$. Thus, $\lk_{\Delta_{t-1}}(v_t)$ is contained in $\st_{\Delta_t}(n-t)$ in $\Delta_t$. Therefore, it is contractible in $\Delta_t$. Applying Lemma \ref{lem:cone} once again yields
\begin{equation}\label{eq:deletion}
    \Delta_{t-1} \simeq  \Delta_t^{\lambda'} \vee \Sigma(\lk_{\Delta_{t-1}}(v_t)).
\end{equation}

\medskip
\noindent\textbf{Step 4: Resolving $\lk_{\Delta_{t - 1}}(v_t)$.}
To resolve $\lk_{\Delta_{t-1}}(v_t)$, consider the vertex $n-t$ in $\lk_{\Delta_{t-1}}(v_t).$ 
It lies in the $(\lambda_1-t)$-th cell of the top row of the diagram supporting $\lk_{\Delta_{t-1}}(v_t)$, and the entries on its right in the top row have been deleted. Moreover, it is the largest vertex in $V(\lk_{\Delta_{t-1}}(v_t)),$ and it only appears in the $(\lambda_1-t)$-th cell of the top row of the diagram. 

One can use an argument similar to that in Step~3 to show that 
\[
\del_{\lk_{\Delta_{t-1}}(v_t)}(n-t) = \Delta_t^{\lambda''},
\]
where $\lambda'' = (\min\{\lambda_1-1, \lambda_2\},\min\{\lambda_1-1, \lambda_3\}, \ldots,\min\{\lambda_1-1, \lambda_r\})$.

On the other hand, by definition, the facets of $\lk_{\lk_{\Delta_{t-1}}(v_t)}(n-t)$ are of the form $F\setminus \{n-t\}$ where $F$ is a facet of $\lk_{\Delta_{t-1}}(v_t)$ containing $n-t$. Therefore,
\[
\lk_{\lk_{\Delta_{t-1}}(v_t)}(n-t) = \Delta_t^{\lambda'''},
\]
where $\lambda''' = (\min\{\lambda_1-t-1, \lambda_1\},\min\{\lambda_1-t-1, \lambda_2\}, \ldots,\min\{\lambda_1-t-1, \lambda_r\})$.

We now claim that $\lk_{\lk_{\Delta_{t-1}}(v_t)}(n-t)$ is contained in $\st_{\mathcal{K}}(n-t-1)$, where $\mathcal{K}=\del_{\lk_{\Delta_{t-1}}(v_t)}(n-t)$.
Let $F$ be a facet of $\lk_{\lk_{\Delta_{t-1}}(v_t)}(n-t)$. 
Then, either $n-t-1 \in F$ or the maximal entry in $F$ is not in the top row of the Young diagram of $\lk_{\lk_{\Delta_{t-1}}(v_t)}(n-t)$ (see Table \ref{tab:link vt} for reference).
In the latter case, observe that $F \cup \{n-t-1\} \in \mathcal{K}= \del_{\lk_{\Delta_{t-1}}(v_t)}(n-t) =\Delta_t^{\lambda''}$. Hence $F\in \st_{\mathcal{K}}(n-t-1)$. This completes the proof of the claim.

Thus, $\lk_{\lk_{\Delta_{t-1}}(v_t)}(n-t)$ is contractible in $\del_{\lk_{\Delta_{t-1}}(v_t)}(n-t)$. Hence, by Lemma \ref{lem:cone}, we have
\begin{equation}\label{eq:link}
    \lk_{\Delta_{t-1}}(v_t) \simeq \Delta_t^{\lambda''} \vee \Sigma(\Delta_t^{\lambda'''}).
\end{equation}
Finally, combining Equations (\ref{eq:deletion}) and (\ref{eq:link}) with $\Delta_t^{\lambda}\simeq \Delta_{t-1}$, we obtain
\[
    \Delta_t^{\lambda} \simeq \Delta_t^{\lambda'} \vee \Sigma(\Delta_t^{\lambda''}) \vee \Sigma^2(\Delta_t^{\lambda'''})
\]
This completes the proof.
\end{proof}

Let $\lambda = (\lambda_1, \lambda_2,  \ldots,  \lambda_r)$ and $\mu = (\mu_1,\mu_2, \ldots, \mu_s)$ be two tuples. We order tuples lexicographically; that is, we write $\mu<\lambda$ if either $s<r$ or $r=s$ and $\mu_i<\lambda_i,$ where $i = \min\{j: \mu_j\neq \lambda_j\}$. We can now show that every $t$-Young complex is either contractible or homotopy equivalent to a wedge of spheres.

\begin{Theorem}\label{thm:homotopytype}
    Let $t$ be a positive integer and $\lambda = (\lambda_1, \lambda_2, \ldots, \lambda_r)$ be a partition.
    Then the simplicial complex $\Delta_t^{\lambda}$ is either contractible or homotopy equivalent to a
wedge of spheres.
\end{Theorem}
\begin{proof}
     We proceed by induction on the partition $\lambda= (\lambda_1,\ldots,\lambda_r)$ in lexicographic order. 
     When $\lambda=(0,0,\ldots,0)$, by our convention $\Delta_t^{\lambda} = \{\emptyset\} = S^{-1}$. 
     If $r = 1,$ then $\Delta_t^{\lambda}$ is a simplex, and hence contractible. Also, if $\lambda=(1,1,\ldots,1)$ is an $r$-tuple, then $\Delta_t^{\lambda}$ is a disjoint union of $r$ points. Thus,
     \[
     \Delta_t^{\lambda} \simeq \mathop{\bigvee}_{r -1} S^0.
     \]
     
     Now assume $r \geq 2$ and that the result holds for all tuples $\mu <\lambda$. If $\lambda_1>\lambda_2$, then $\Delta_t^{\lambda}$ is a cone with apex $(r-1)t+\lambda_1$, and is therefore contractible. Therefore, we may assume that $\lambda_1= \lambda_2$. 
     
     When $\lambda_1\leq t,$ by Lemma~\ref{lem:wedgeproduct1} we have
    \[
    \Delta_t^{\lambda}\simeq \Delta_t^{\lambda'} \vee \Sigma(\Delta_t^{\lambda''}),
    \] 
    where $\lambda' = (\lambda_2, \lambda_3,  \ldots,  \lambda_r),$ and $\lambda'' = (\min\{\lambda_1-1, \lambda_2\},\min\{\lambda_1-1, \lambda_3\}, \ldots,\min\{\lambda_1-1, \lambda_r\})$. Since $\lambda'<\lambda$ and $\lambda''<\lambda$, the induction hypothesis applies. Thus, both $\Delta_t^{\lambda'}$ and $\Delta_t^{\lambda''}$ are either contractible or homotopy equivalent to a
wedge of spheres. Because the suspension of a contractible space is contractible, and the suspension of a wedge of spheres is again a wedge of spheres, it follows that $\Delta_t^{\lambda}$ is either contractible or homotopy equivalent to a
wedge of spheres.

On the other hand, when $\lambda_1 > t$, under the notation of Lemma~\ref{lem:wedgeproduct2}, we have that $\lambda'<\lambda,$ $\lambda''<\lambda,$ and $\lambda'''<\lambda$. By the same inductive argument as above, the claim follows in this case as well.
\end{proof}

We now investigate the vertex-decomposability of $t$-Young complexes.
We first show that every $1$-Young complex is vertex-decomposable and then characterize all vertex-decomposable $t$-Young complexes for $t \ge 2$.

\begin{Theorem}\label{thm:t=1vd}
Let $\lambda = (\lambda_1, \lambda_2, \ldots, \lambda_r)$ be a partition. Then, $\Delta_1^{\lambda}$ is vertex-decomposable.
\end{Theorem}

\begin{proof}
     

We proceed by induction on the well-order of partitions introduced above. Thus, for a given partition $\lambda=(\lambda_1,\ldots,\lambda_r)$, we assume that the assertion holds for every partition $\mu<\lambda$.

If $r=1$, then $\Delta^\lambda_1$ is a simplex and hence is vertex-decomposable.

We next consider the case $\lambda=(1,\ldots,1)$. In this case, $\Delta^\lambda_1$ is a zero-dimensional simplicial complex consisting of $r$ isolated vertices. Therefore, $\Delta^\lambda_1$ is vertex-decomposable.

We may now assume that $r\geq2$ and that $\lambda\neq(1,\ldots,1)$. Suppose first that $\lambda_1>\lambda_2$. Then the top two rows of the Young diagram corresponding to $\lambda$ are shown in Table~\ref{tab: top two rows vertex decomposibility}.
    
\begin{table}[h]
\begin{tabular}{ |c|c|c|c|c|c| } 
 \hline
 $(r-1)+ 1$ & $\cdots$ & $(r-1)+ \lambda_2$ & \cellcolor{lightgray}$(r-1)+ (\lambda_2+1)$ & \cellcolor{lightgray}$\cdots$ & \cellcolor{lightgray}$(r-1)+ \lambda_1$ \\
\hline
 $(r-2)+ 1$ & $\cdots$ & $(r-2)+ \lambda_2$ & \multicolumn{3}{c}{} \\
 \cline{1-3}
\end{tabular}
\caption{Top two rows of Young diagram when $\lambda_1>\lambda_2$.}
\label{tab: top two rows vertex decomposibility}
\end{table}
Therefore, we have 
\[\Delta_1^{\lambda} =\Delta_1^{\mu} * \big\langle\{(r-1) + \lambda_2+1,\ldots,(r-1)+ \lambda_1\}\big\rangle,\]
where $\mu = (\lambda_2, \lambda_2, \lambda_3, \ldots, \lambda_r)$.
The partitions $\mu$ and $\lambda$ have the same number of components,
and the first component of $\mu$ is
$\lambda_2$, which is less than $\lambda_1$. Hence $\mu<\lambda$. By the induction hypothesis, $\Delta^\mu_1$ is vertex-decomposable. Since the join of
a vertex-decomposable complex with a simplex is vertex-decomposable \cite[Proposition~2.4]{PB}, it follows that
$\Delta^\lambda_1$ is vertex-decomposable.

We may therefore assume that $\lambda_1=\lambda_2$. Since $\lambda\neq(1,\ldots,1)$, we have $\lambda_1\geq2$. In this case, the top two rows of the Young diagram of $\lambda$ are shown in Table~\ref{tab: top two rows vertex decomposibility 2}.

\begin{table}[h]
\begin{tabular}{ |c|c|c|c |c| } 
 \hline
 $(r-1)+ 1$ & $(r-1)+ 2$ & $\cdots$ & $(r-1)+ (\lambda_1-1)$ &   $(r-1) + \lambda_1$\\
 \hline
 $(r-1) $ & $(r-1)+ 1$ & $\cdots$ & $(r-1)+ (\lambda_1-2)$ &   $(r-1)+ (\lambda_1-1)$ \\
  \hline
\end{tabular}
\caption{Top two rows of Young diagram when $\lambda_1=\lambda_2$.}
\label{tab: top two rows vertex decomposibility 2}
\end{table}

Set $v = (r-1)+ \lambda_1$. The vertex $v$ appears only in the last cell of the top row.
Following the argument of Step~3 in the proof of Lemma~\ref{lem:wedgeproduct2}, we obtain
\[
\del_{\Delta_1^{\lambda}}(v) =\Delta_1^{\lambda'},\quad \text{where} \ \lambda' = (\lambda_2, \lambda_3, \ldots, \lambda_r).
\]
Indeed, it is clear that $\Delta_1^{\lambda'} \subseteq \del_{\Delta_1^{\lambda}}(v)$. Conversely, let $F$ be any face of $\del_{\Delta_1^{\lambda}}(v)$ containing one of the top-row vertices, say $(r-1)+j \in \MS_1$ for some $1\leq j\leq \lambda_1 -1$. The entries of $F$ that are greater than $(r-1)+j$ correspond to the vertices to the right in the top row.
We may reinterpret this vertex $(r-1)+j$ as an entry in the $(j+1)$-th cell of the second row from the top, and similarly, if $F$ contains the vertex in the $(j+1)$-th cell of the top row, we can view it as corresponding to the $(j+2)$-th cell of the second row, and so on.  
Since there is no vertex in the last cell of the top-row (as $v$ has been deleted), $F$ can be realized as a face in the simplicial complex obtained after deleting the top row. Thus, $F\in \Delta_1^{\lambda'}$.

Since $\lambda_1=\lambda_2$, we have $\dim(\Delta_1^{\lambda'})=\dim(\Delta_1^{\lambda}) = \lambda_1-1$, so $v$ is a shedding vertex of $\Delta_1^{\lambda}$. Moreover, since $\lambda'<\lambda$, by the induction hypothesis $\del_{\Delta_1^{\lambda}}(v)$ is vertex-decomposable.
We now describe the link. Following again the argument in the proof of Lemma~\ref{lem:wedgeproduct2}, we have  
\[
\lk_{\Delta_1^{\lambda}}(v) =\Delta_1^{\lambda''}, \quad \text{where} \ \lambda'' = (\lambda_1-1, \lambda_2-1, \min\{\lambda_3,\lambda_1-1\},  \ldots, \min\{\lambda_r,\lambda_1-1\}).
\]
Since $\lambda'' <\lambda$, by the induction hypothesis $\lk_{\Delta_1^{\lambda}}(v)$ is vertex-decomposable.

Thus, $v$ is a shedding vertex of $\Delta^\lambda_1$, and both $\del_{\Delta^\lambda_1}(v)$ and
$\lk_{\Delta^\lambda_1}(v)$ are vertex-decomposable. Therefore,
$\Delta^\lambda_1$ is vertex-decomposable.
\end{proof}

\begin{Theorem}\label{thm:cmdiagram}
Let $t\geq 2$ and $\lambda = (\lambda_1, \lambda_2, \ldots, \lambda_r)$ be a partition. Then the following are equivalent:
\begin{enumerate}
    \item \label{thm:cmdiagram:vd} $\Delta_t^{\lambda}$ is vertex-decomposable.
    \item \label{thm:cmdiagram:shellable} $\Delta_t^{\lambda}$ is shellable.
    \item \label{thm:cmdiagram:cm} $\Delta_t^{\lambda}$ is Cohen-Macaulay.
    \item \label{thm:cmdiagram:t} Either $r=1$ or $r\geq 2$ and $\lambda_2 \leq t$.
\end{enumerate}
\end{Theorem}

\begin{proof}
The implications $\eqref{thm:cmdiagram:vd} \implies \eqref{thm:cmdiagram:shellable}$ and $\eqref{thm:cmdiagram:shellable} \implies \eqref{thm:cmdiagram:cm}$ always hold in general, as established in \cite[Theorem 11.3]{BW97part2} and \cite[Theorem 5.1.13]{BH}, respectively.
 
$\eqref{thm:cmdiagram:cm} \implies \eqref{thm:cmdiagram:t}$
We first claim that if $\lambda_1= \lambda_2 =t+1,$ then $\Delta_t^{\lambda}$ is not Cohen-Macaulay. 
By Lemma~\ref{lem:wedgeproduct2}, 
 \[
    \Delta_t^{\lambda}\simeq \Delta_t^{\lambda'} \vee \Sigma(\Delta_t^{\lambda''}) \vee \Sigma^2(\Delta_t^{\lambda'''}),
    \] 
where $\lambda', \lambda''$ and $\lambda'''$ are as defined in Lemma~\ref{lem:wedgeproduct2}. We are only interested in 
\[
\lambda'''=(\min\{\lambda_1-t-1, \lambda_1\},\min\{\lambda_1-t-1, \lambda_2\}, \ldots,\min\{\lambda_1-t-1, \lambda_r\}).
\]
Since $\lambda_1 = t+1$ and $\lambda_1\geq \lambda_i$ for all $1\leq i \leq r,$ we get that $\lambda''' = (0,0,\ldots,0)$. 
Thus, $\Delta_t^{\lambda'''} = \{\emptyset\}$, hence $\Sigma^2(\Delta_t^{\lambda'''}) = S^1.$ In particular, $\widetilde{H}_1(\Delta_t^{\lambda}) \neq 0.$ 
Since $\Delta_t^{\lambda}$ is a simplicial complex of dimension $\lambda_1 - 1 = t$ with $t \geq 2$, it follows that it has a nontrivial homology in dimensions strictly less than~$t$. Consequently, $\Delta_t^{\lambda}$ is not Cohen--Macaulay.

Now suppose that $\lambda_2 \geq t+1$ and $\lambda_1\neq t+1$.
Since $\lambda_1\geq \lambda_2$ and $\lambda_1\neq t+1$, the top row of the Young diagram $\lambda$ is
\begin{center}
\begin{tabular}{ |c|c|c|c |c|c| } 
 \hline
 $(r-1) t+ 1$ & $\cdots$ & $(r-1) t+ (t+1)$ &  \cellcolor{lightgray} $(r-1) t+ (t+2)$ &\cellcolor{lightgray} $\cdots$ & \cellcolor{lightgray} $(r-1) t+ \lambda_1$\\
 \hline
\end{tabular}.
\end{center}
Consider the face $F = \{(r-1) t+ (t+2),\ldots,(r-1) t+ \lambda_1\}$ in $\Delta_t^{\lambda}$. Using an argument similar to that in Lemma~\ref{lem:wedgeproduct2}, any maximal admissible choice from the first $t+1$ columns may appear in a face of $\Delta_t^{\lambda}$ containing $F$. Therefore,
 \[\lk_{\Delta_t^{\lambda}}(F) = \Delta_t^{\mu},\] 
 where $\mu = (\min\{t+1, \lambda_1\},\min\{t+1, \lambda_2\}, \ldots,\min\{t+1, \lambda_r\})$. Since $\lambda_1\geq \lambda_2\geq t+1$,
 we have $\mu = (t+1, t+1, \min\{t+1, \lambda_3\} \ldots,\min\{t+1, \lambda_r\})$. By the previous claim, $\Delta_t^{\mu}$ is not Cohen-Macaulay. Since the link of a Cohen–Macaulay simplicial complex must itself be Cohen–Macaulay, we conclude that $\Delta_t^{\lambda}$ is not Cohen-Macaulay. This completes the proof.

 $\eqref{thm:cmdiagram:t} \implies \eqref{thm:cmdiagram:vd}$ Suppose that $\lambda_2 \leq t.$ This means that each entry in the Young diagram $\lambda$ appears in exactly one cell of $\lambda$.  We proceed by induction on $\lambda= (\lambda_1,\ldots,\lambda_r)$, ordered lexicographically. If $r=1$, then $\Delta_t^{\lambda}$ is a simplex and hence vertex-decomposable. Now assume $r\geq 2$ and that the claim holds for all partitions $\mu<\lambda$. 

Set $v:=\lambda_r \in V(\Delta_t^{\lambda})$, the entry in the last cell of the bottom (i.e. $r$-th) row. When $v\geq 2$, the other entries of the bottom row are $1, 2, \ldots, v-1$ from left to right. We first show that $v$ is a shedding vertex of $\Delta_t^{\lambda}$. 
 Let $F$ be any facet of $\Delta_t^{\lambda}$ containing $v$.
By definition of $\Delta_t^\lambda$, each facet selects exactly one entry from each column, and the chosen entries increase strictly from left to right. In particular $F$ contains the initial block $\{1,2,\dots,v-1\}$ to the left of $v$  together with a maximal admissible choice from the columns $\lambda_r+1$ onward.
Since $\lambda_{r-1}\geq \lambda_r$, the cell immediately above the position of $v$ exists in the Young diagram. The entry in this cell is $\lambda_r+t$, which we denote by $v'$. Because $v$ is the last entry of the bottom row, the set $(F\setminus \{v\}) \cup \{v'\}$ is also a facet of $\Delta_t^{\lambda}$. Thus, $v$ is a shedding vertex of $\Delta_t^{\lambda}$.

\begin{table}[h]
\centering
\begin{tabular}{ |c|c|c|c | c |c |  }
\hline
 $(r-1) t+ 1$ & $\cdots$ & $(r-1) t+ v$ &  \cellcolor{lightgray} $(r-1) t+ (v+1)$ &\cellcolor{lightgray} $\cdots$ & \cellcolor{lightgray} $(r-1) t+ \lambda_1$ \\
\hline
$(r-2) t+ 1$ &  $\cdots$ &  $(r-2) t+ v$ & \cellcolor{lightgray}$(r-2) t+ (v+1)$ &\cellcolor{lightgray} $\cdots$ & \multicolumn{1}{c}{}  \\
\cline{1-5}
 $\vdots$ & $\cdots$ &  $\vdots$ &\cellcolor{lightgray} $\vdots$ & \cellcolor{lightgray}$\cdots$ & \multicolumn{1}{c}{} \\
\cline{1-5} 
1 & $\cdots$ & $v$ &  \multicolumn{3}{c}{}  \\
\cline{1-3}
\end{tabular}
\caption{$\Delta_t^\tau$ is isomorphic to the subcomplex with faces defined using entries from the shaded region.}
\end{table}

As observed above, every facet $F$ containing $v$ includes $\{1,2,\dots,v\}$ and a maximal admissible choice from the columns $\lambda_r+1$ onward. Therefore
\[
\lk_{\Delta_t^\lambda}(v) = \Delta_1 * \big\langle\{1,2,\dots,v-1\}\big\rangle,
\]
where $\Delta_1$ is the subcomplex of $\Delta_t^\lambda$ consisting of all admissible faces supported in the columns indexed $\lambda_r+1,\lambda_r+2,\dots,\lambda_1$. Removing the first $\lambda_r$ columns from $\lambda$ yields a diagram whose $j$-th row has length $\lambda_j-\lambda_r$ for $1\leq j\leq r-1$. Since the properties of a $t$-Young complex depend on the structure of partition and how the entries are labeled, i.e., increase by $t$ while going up in a column and increase by one while going right in a row, we get that $\Delta_1$ is isomorphic to $\Delta_t^{\tau},$ where $\tau = (\lambda_1-\lambda_r, \ldots, \lambda_{r-1}-\lambda_r).$ 
Since $\tau<\lambda$, the induction hypothesis implies that \(\Delta_t^{\tau}\) is vertex-decomposable. Because the join of a vertex-decomposable complex with a simplex is vertex-decomposable, it follows that $\lk_{\Delta_t^\lambda}(v)$ is vertex-decomposable.

Deleting $v$ amounts to shortening the last row by one cell while leaving the other row lengths unchanged. Therefore, the deletion corresponds to the partition $\mu := (\lambda_1,\lambda_2,\ldots,\lambda_{r-1},\lambda_r-1)$. More precisely, $\del_{\Delta_t^\lambda}(v)=\Delta_t^\mu$ when $\lambda_r\geq 2$, and  $\del_{\Delta_t^\lambda}(v)$ is isomorphic to $\Delta_t^\mu$ when $\lambda_r=1$. Because $\mu<\lambda$, the induction hypothesis guarantees that $\del_{\Delta_t^\lambda}(v)$ is vertex-decomposable. Hence, the proof is complete.
\end{proof}

\section{From $t$-Young complexes to $t$-path ideals of path graphs}\label{Section: Specialization to t-path} 

In this section, we show how the topological results established in the previous section can be specialized to a combinatorial–algebraic setting, in particular to the case of $t$-path ideals of path graphs and their squarefree powers. 

Let us begin with establishing in the following proposition a connection between $t$-Young complexes and the Alexander duals of the Stanley--Reisner complexes associated with $I_{n,t}^{[k]}$. 

\begin{Proposition}\label{lem:dualfacets}
Under the assumptions of Notation~\ref{Notation}, we have
\[
\Delta_{n,t}^{[k]}
=
\begin{cases}
\emptyset, & \text{if } n<kt,\\[0.3em]
\{\emptyset\}, & \text{if } n=kt,\\[0.3em]
\Delta_t^\lambda, & \text{if } n>kt,
\end{cases}
\]

where
$\lambda=(\lambda_1,\ldots,\lambda_{k+1})$ is the partition defined
by $\lambda_i=n-kt$ for every $1\leq i\leq k+1$.
\end{Proposition}

\begin{proof}
Suppose first that $n<kt$. Then $I_{n,t}^{[k]}=(0)$, and hence $\Sigma_{n,t}^{[k]}$ is the full simplex on $[n]$. Therefore, $\Delta_{n,t}^{[k]}=\emptyset.$

Suppose now that $n=kt$. In this case, $I_{n,t}^{[k]}=(x_1x_2\cdots x_n)$, so it easily follows that $\Delta_{n,t}^{[k]}
    =
    \{[n]\setminus[n]\}
    =
    \{\emptyset\}.$ 

Assume now that $n > kt$. We show that $\Delta_{n,t}^{[k]} = \Delta_t^{\lambda}$, where $\lambda = (\lambda_1, \lambda_2, \ldots, \lambda_{k+1})$ with $\lambda_i = n - kt$, for all $1\leq i\leq k+1$.  
The corresponding Young diagram has $k+1$ rows, each of length $n - kt$, and is represented as in Table \ref{tab: Young diagram for t-path}.
\begin{table}[h]
\centering
\begin{tabular}{ |c|c|c|c|c| } 
 \hline
 $kt+1$ & $kt+2$ & $\cdots$ & $n-1$ & $n$ \\
 \hline
 $(k-1)t+1$ & $(k-1)t+2$ & $\cdots$ & $n-t-1$ & $n-t$ \\ 
 \hline
 $\vdots$ & $\vdots$ & $\cdots$ & $\vdots$ & $\vdots$ \\ 
 \hline
 $1$ & $2$ & $\cdots$ & $n-kt-1$ & $n-kt$ \\ 
 \hline
\end{tabular}
\caption{Young diagram of $\Delta_{n,t}^{[k]}$.}
\label{tab: Young diagram for t-path}
\end{table}

By the definition of $I_{n,t}^{[k]}$ and $\Sigma_{n,t}^{[k]}$, an element of $\MN(\Sigma_{n,t}^{[k]})$ can be written as
\[
\bigcup_{i=1}^k \{\,b_i, b_i + 1, \dots, b_i + t - 1\,\},
\quad \text{where } b_i + t - 1 < b_{i+1} \text{ for all } 1 \le i \le k - 1.
\]
Let $F$ be a facet of $\Delta_{n,t}^{[k]}$. It is known that 
$\mathcal{F}(\Delta_{n,t}^{[k]}) = \big\{ [n] \setminus F : F \in \MN(\Sigma_{n,t}^{[k]}) \big\}$ (see \cite[page 17]{HHmonomialideals}). 
Hence, writing $F=\{a_1<a_2<\cdots<a_{n-kt}\},$ we have $a_1\equiv1\pmod t$ and $a_{i+1}\equiv a_i+1\pmod t$ for every $1\leq i\leq n-kt-1$. Equivalently, each $a_j$ lies in the $j$-th column of the Young diagram in Table~\ref{tab: Young diagram for t-path}.


Therefore, $F$ is a facet of $\Delta_{n,t}^{[k]}$ if and only if $F$ is a facet of $\Delta_t^{\lambda}$. Hence $\Delta_{n,t}^{[k]} = \Delta_t^{\lambda}$, as desired. 
\end{proof}

Note that for $t=1$, the simplicial complex $\Delta_{n,1}^{[k]}$ is the $(n-k-1)$-skeleton of the $(n-1)$-simplex, and therefore the ideal $I_{n,1}^{[k]}$ is the \emph{squarefree Veronese ideal} of degree $k$. Moreover, as a consequence of the above proposition, we compute the Helly number of Alexander dual of all $t$-Young complexes. Recall that the \emph{Helly number} of a simplicial complex $\Delta$, denoted by $h(\Delta)$, is the maximum cardinality of a minimal non-face of $\Delta$. Equivalently, $h(\Delta)$ is the maximum degree of a minimal monomial generator of the Stanley--Reisner ideal $I_\Delta$.

\begin{Corollary}\label{cor:helly}
    Let $t$ be a positive integer and $\lambda = (\lambda_1, \lambda_2, \ldots, \lambda_r)$ be a partition with $r\geq 2$. Set $n = (r-1)t+\lambda_1$. Let $\Sigma_{t}^{\lambda}$ denote the Alexander dual of $\Delta^{\lambda}_t$ on vertex set $[n]$. Then, the Stanley-Reisner ideal of $\Sigma_{t}^{\lambda}$ is generated by a subset of generators of $I_{n,t}^{[r-1]}$. In particular, the Helly number of $\Sigma_{t}^{\lambda}$ is
    \[
    h({\Sigma^{\lambda}_t}) = (r-1)t.
    \]
\end{Corollary}

\begin{proof}
  We claim that every facet of $\Delta^{\lambda}_t$ is also a facet of $\Delta_{n,t}^{[r-1]}$. 
  By Proposition~\ref{lem:dualfacets}, $\Delta_{n,t}^{[r-1]} = \Delta_t^{\mu}$ where $\mu = (\mu_1,\ldots,\mu_r)$ with $\mu_i = n - (r-1)t=\lambda_1$, for all $1\leq i\leq r$.
  Note that $\lambda\leq \mu$, i.e., every cell in $\lambda$ also appears in $\mu$ at exactly same position. Moreover, the entry of each such cell is same in both $\lambda$ and $\mu$. 
  Since for a facet $F= \{a_1< \dots< a_{\lambda_1}\}$ of $\Delta^{\lambda}_t$, we have that each $a_j$ lies in the $j$-th column of the Young diagram $\lambda$, for all $1 \le j \le \lambda_1$. By above observations, $F$ is also a face of $\Delta_{n,t}^{[r-1]}$.
  Hence the claim follows as $\dim (\Delta^{\lambda}_t) = \dim(\Delta_{n,t}^{[r-1]}) =\lambda_1-1$. Therefore, $\mathcal{F}(\Delta^{\lambda}_t) \subseteq \mathcal{F}(\Delta_{n,t}^{[r-1]})$. 
  
  Taking complements in $[n]$ and using the definition of the
Alexander dual, we obtain
\[
\mathcal{N}(\Sigma_t^\lambda)
=
\bigl\{[n]\setminus F:
F\in\mathcal{F}(\Delta_t^\lambda)\bigr\}\\
\subseteq
\bigl\{[n]\setminus F:
F\in\mathcal{F}(\Delta_{n,t}^{[r-1]})\bigr\}\\
=
\bigl\{\operatorname{supp}(u):
u\in G(I_{n,t}^{[r-1]})\bigr\}.
\]

It follows that every minimal nonface of $\Sigma_t^\lambda$ is the support of a minimal generator of $I_{n,t}^{[r-1]}$. Therefore, the Stanley--Reisner ideal of $\Sigma_t^\lambda$ is generated by a subset of the minimal generators of $I_{n,t}^{[r-1]}$. In particular, every minimal generator of $I_{n,t}^{[r-1]}$ is the product of $r-1$ pairwise disjoint $t$-paths and hence has degree $(r-1)t$. Thus every minimal non-face of $\Sigma_t^\lambda$ has cardinality $(r-1)t$, so we have that $h(\Sigma_t^\lambda)=(r-1)t,$ as desired.

\end{proof}

Using the structural characterization of $\Delta_{n,t}^{[k]}$ established in Proposition~\ref{lem:dualfacets}, we can now express the homotopy type of $\Delta_{n,t}^{[k]}$ in terms of lower-dimensional complexes.

\begin{Corollary}\label{cor:formulahomotopy}
Let $n, k\in \NN$ and $t$ be a positive integer such that $n = kt+l$ for some $l\geq 1.$
\begin{enumerate}
    \item if $l\leq t$, then 
    \[
    \Delta_{n,t}^{[k]} \simeq \Delta_{n-t,t}^{[k-1]}  \vee \Sigma(\Delta_{n-t-1,t}^{[k-1]}).
    \]
    \item if $l > t$, then 
    \[
\Delta_{n,t}^{[k]} 
\simeq 
\Delta_{n-t,t}^{[k-1]} 
\vee 
\Sigma(\Delta_{n-t-1,t}^{[k-1]}) 
\vee 
\Sigma^2(\Delta_{n-t-1,t}^{[k]}),
\]
where $\Delta_{kt,t}^{[k]}= \{\emptyset\}$.
\end{enumerate}
\end{Corollary}

\begin{proof}
    The proof follows from Lemma~\ref{lem:wedgeproduct1}, Lemma~\ref{lem:wedgeproduct2} and Proposition~\ref{lem:dualfacets}
\end{proof}

In \cite{KNQ24squarefree}, all squarefree powers of $t$-path ideals of path graphs having a linear resolution are characterized in terms of the matching and restricted matching numbers. In the following remark, we recall that result and observe that it follows directly from the results obtained so far, in particular from Theorem \ref{thm:cmdiagram} and Proposition~\ref{lem:dualfacets}.

\begin{Remark}\label{Remark: linear resolution}
Recalling Notation~\ref{Notation}, let $\Gamma_{n,t}$ denote the facet complex of $I_{n,t}$ and let $\Delta_{n,t}^{[k]}$ be the Alexander dual of the Stanley–Reisner complex of $I_{n,t}^{[k]}$.

From \cite[Theorems~2.7 and~3.3]{KNQ24squarefree}, the following statements are equivalent for $t \geq 2$:
\begin{enumerate}
\item[(i)] $I_{n,t}^{[k]}$ has linear quotients;
\item[(ii)] $I_{n,t}^{[k]}$ has a linear resolution;
\item[(iii)] $k = \nu_0(\Gamma_{n,t})$ or $k = \nu(\Gamma_{n,t})$.
\end{enumerate}
where $\nu_0(\Gamma_{n,t})$ and $\nu(\Gamma_{n,t})$ denote the restricted matching number and the matching number of $\Gamma_{n,t}$, respectively.
Moreover, by \cite[Remark~4.2]{KNQ24squarefree},
    \[
        \nu_0(\Gamma_{n,t}) = \left\lfloor \frac{n - 1}{t} \right\rfloor
        \quad \text{and} \quad
        \nu(\Gamma_{n,t}) = \left\lfloor \frac{n}{t} \right\rfloor.
    \]
   Hence we have that
    \begin{itemize}
        \item $k = \nu(\Gamma_{n,t}) \ \text{if and only if} \ kt\leq n\leq kt+t-1$;
        \item $k = \nu_0(\Gamma_{n,t})     \ \text{if and only if} \
        kt+1 \leq n \leq kt + t.$
    \end{itemize}
    
 It is known from \cite[Theorem~3]{ER98linearresolution} that a squarefree monomial ideal has a linear resolution if and only if the Alexander dual of its Stanley–Reisner complex is Cohen–Macaulay. Hence, $I_{n,t}^{[k]}$ has a linear resolution if and only if $\Delta_{n,t}^{[k]}$ is Cohen–Macaulay.

Similarly, a squarefree monomial ideal has linear quotients if and only if its Alexander dual complex is shellable (see \cite[Proposition~8.2.5]{HHmonomialideals}). Therefore, $I_{n,t}^{[k]}$ has linear quotients if and only if $\Delta_{n,t}^{[k]}$ is shellable.

 Therefore, the equivalence of (i), (ii), and (iii) above follow directly by combining \cite[Theorem~3]{ER98linearresolution},  \cite[Proposition~8.2.5]{HHmonomialideals}, Theorem \ref{thm:cmdiagram} and Proposition~\ref{lem:dualfacets}.
\end{Remark}

Taking into account the results discussed in the previous remark, we obtain the following corollary as a consequence of Theorem \ref{thm:t=1vd} and Theorem \ref{thm:cmdiagram}.

\begin{Corollary}\label{cor:linearity}
    Let $\Sigma^{\lambda}_t$ be as in Corollary~\ref{cor:helly} and $I_{\Sigma^{\lambda}_t}$ be its Stanley-Reisner ideal. Then
    \begin{enumerate}
        \item If $t=1$, then $I_{\Sigma^{\lambda}_t}$ has linear quotients.
        \item If $t\geq 2$, then $I_{\Sigma^{\lambda}_t}$ has linear quotients if and only if it has a linear resolution if and only if $\lambda_2\leq t$.
    \end{enumerate}
\end{Corollary}

Now we turn to the homotopy type of the $t$-Young complexes $\Delta_{n,t}^{[k]}$.  We recall from Theorem \ref{thm:homotopytype} that $\Delta_{n,t}^{[k]}$ is either contractible or homotopy equivalent to a wedge of spheres.
In the following, we determine the homotopy type of $\Delta_{n,t}^{[k]}$ explicitly for all $k,t\geq 1$ and $n>kt$. Set $s=n-kt$. \cref{cor:formulahomotopy} will be the key for determining the homotopy type of $\Delta_{n,t}^{[k]}$.

For $j,u\geq 0$, define
\[
p_{j,u}(z)=\sum_{i\geq -1}\dim_K \widetilde{H}_i\big(\Delta_{jt+u,t}^{[j]}\big)\,z^{i+1},
\]
and set $p_{j,u}(z)=0$ for $u<0$. Since $\Delta_{jt+u,t}^{[j]}$ is a finite dimensional simplicial complex, $p_{j,u}(z)$ is a polynomial in $z$. Note that $p_{j,0}(z)=1$ for all $j\geq 0$, because $\Delta_{jt,t}^{[j]}=\{\emptyset\}$ and $\widetilde{H}_{-1}(\{\emptyset\})=K$, and that $p_{0,u}(z)=0$ for all $u\geq 1$, because $\Delta_{u,t}^{[0]}$ is a simplex.

Since suspension shifts reduced homology by one degree, \cref{cor:formulahomotopy} yields
\begin{equation}\label{eq:homotopy}
p_{j,u}(z)=p_{j-1,u}(z)+z\,p_{j-1,u-1}(z)+z^2\,p_{j,u-t-1}(z)\qforall j,u\geq 1\,.
\end{equation}
Observe that for $1\leq u\leq t$ the last term vanishes, since its second index is negative.

For $j\geq 0$, define the formal power series in variable $x$ with coefficients in $\ZZ[z]$
\[
F_j(x,z)=\sum_{u\geq 0}p_{j,u}(z)\,x^u\in\ZZ[z][[x]]\,.
\]
Then $F_0(x,z)=1$, and for $j\geq 1$ we have
\begin{align*}
\sum_{u\geq 1}p_{j-1,u}(z)\,x^u &= F_{j-1}(x,z)-1\,,\\
\sum_{u\geq 1}z\,p_{j-1,u-1}(z)\,x^u &= zx\,F_{j-1}(x,z)\,,\\
\sum_{u\geq 1}z^2\,p_{j,u-t-1}(z)\,x^u &= z^2x^{t+1}\,F_j(x,z)\,,
\end{align*}
where the last identity uses $p_{j,v}(z)=0$ for $v<0$. Hence, by \eqref{eq:homotopy}, for $j\geq 1$,
\begin{align*}
F_j(x,z)&=1+\sum_{u\geq 1}p_{j,u}(z)\,x^u\\
&=1+\sum_{u\geq 1}\big(p_{j-1,u}(z)+z\,p_{j-1,u-1}(z)+z^2\,p_{j,u-t-1}(z)\big)x^u\\
&=1+\big(F_{j-1}(x,z)-1\big)+zx\,F_{j-1}(x,z)+z^2x^{t+1}F_j(x,z)\\
&=(1+xz)\,F_{j-1}(x,z)+z^2x^{t+1}F_j(x,z)\,,
\end{align*}
that is,
\[
(1-z^2x^{t+1})\,F_j(x,z)=(1+xz)\,F_{j-1}(x,z)\,.
\]

The element $1-z^2x^{t+1}$ is invertible in $\ZZ[z][[x]]$, with inverse $\sum_{q\geq 0}z^{2q}x^{(t+1)q}$. Multiplying by this inverse we obtain
\[
F_j(x,z)=\frac{1+xz}{1-z^2x^{t+1}}\,F_{j-1}(x,z)\qforall j\geq 1\,.
\]
The factor on the right hand side is independent of $j$. Starting from $F_0=1$ and applying this identity successively for $j=1,\ldots,k$, we get
\[
F_k(x,z)=\Big(\frac{1+xz}{1-z^2x^{t+1}}\Big)^k=\frac{(1+xz)^k}{(1-z^2x^{t+1})^k}\,.
\]
Since
\[
(1+xz)^k=\sum_{b=0}^{k}\binom{k}{b}z^bx^b\qand (1-z^2x^{t+1})^{-k}=\sum_{q\geq 0}\binom{k+q-1}{q}z^{2q}x^{(t+1)q}\,,
\]
we obtain
\[
F_k(x,z)=\sum_{b=0}^{k}\sum_{q\geq 0}\binom{k}{b}\binom{k+q-1}{q}\,z^{b+2q}\,x^{b+(t+1)q}\,.
\]

Now $p_{k,s}(z)$ is the coefficient of $x^s$ in $F_k(x,z)$. A term of this double sum contributes to the coefficient of $x^s$ if and only if $b+(t+1)q=s$; so for each fixed $q$ the only possible value is $b=s-(t+1)q$, and the condition $b\geq 0$ forces $0\leq q\leq\lfloor s/(t+1)\rfloor$. (Terms with $b>k$ vanish, since then $\binom{k}{b}=0$.) For such terms, $b+2q=s-(t+1)q+2q=s-(t-1)q$. Hence
\[
p_{k,s}(z)=\sum_{q=0}^{\lfloor s/(t+1)\rfloor}\binom{k+q-1}{q}\binom{k}{s-(t+1)q}\,z^{s-(t-1)q}\,.
\]
Finally, we conclude using \cref{thm:homotopytype}:
\begin{Theorem}\label{thm:homotopyformula}
    Let $k,t\geq 1$ and $n>kt$ be integers, and set $s=n-kt$. Then
    \[
    \Delta_{n,t}^{[k]} \simeq \bigvee_{q=0}^{\lfloor s/(t+1)\rfloor} \bigvee_{{k+q-1 \choose q}{k \choose s-(t+1)q}}S^{s-(t-1)q-1}\,,
    \]
    where the empty wedge is a point.
\end{Theorem}

\begin{Corollary}\label{prop:homology}
Let $k,t\geq 1$ and $n>kt$ be integers, and set $s=n-kt$. If $s\leq t$, then
    \[
    \Delta_{n,t}^{[k]} \simeq \bigvee_{{k \choose s}} S^{s-1}, 
    \]
    where the right hand side is a point for $k<s$.
\end{Corollary}

We now illustrate the above theorem in the following example:

\begin{Example}\label{eg:homotopy}
    \rm
    The Young diagram of the partition $(3,3,3,3)$ corresponding to $n=9$, $k=3$ and $t =2$ is shown in Table \ref{tab:example}. We have $s=3>2=t$ and $\lfloor s/(t+1)\rfloor =1$. Therefore, by \cref{thm:homotopyformula}
    \begin{align*}
    \Delta_{9,2}^{[3]}& \quad \simeq \quad  \Big( \bigvee_{\binom{3+0-1}{0}\binom{3}{3-(2+1)0}}\,S^{3-(2-1)0-1} \Big) \quad \bigvee \quad
    \Big(
    \bigvee_{\binom{3+1-1}{1}\binom{3}{3-(2+1)1}}\,S^{3-(2-1)1-1} \Big)\\
    & \quad = \quad S^2 \, \bigvee \, \big(\bigvee_3 S^1 \big )\,.
    \end{align*}
\end{Example}

\begin{table}[ht]
\centering
\begin{ytableau}
7 & 8 & 9 \\
5 & 6 & 7\\
3 & 4 & 5 \\ 
1 & 2 & 3 \\ 
\end{ytableau}
\caption{Young diagram filled by $\MS_1$, $\MS_2$, $\MS_3$ and $\MS_4$.}
\label{tab:example}
\end{table}

\section{Algebraic invariants of squarefree powers of \texorpdfstring{$t$}{t}-path ideals of path graphs}\label{Section: Projective dimension}

In this section, we provide a combinatorial description of the projective dimension (Theorem \ref{Theorem: Projective dimension}) and Krull dimension (Theorem \ref{Theorem: Krull dimension}) of the $k$-th squarefree power of the $t$-path ideal of a path graph $P_n$, in terms of the parameters $n$, $t$, and $k$.

Here below, we state the formula for projective dimension.

     \begin{Theorem}\label{Theorem: Projective dimension}
    Under the Notation \ref{Notation}, for any $1\leq k \leq \nu(\Gamma_{n,t})$, we have
        $$ 
        \mathrm{pd}\left(R/I_{n,t}^{[k]}\right)= \begin{cases}
         n-kt+1 & \text{ if }  kt \leq n \leq k(t+1); \\[0.6em]
         \displaystyle\frac{2(n-d)}{t+1}-k+1 & \text{ if } n \equiv d \mod(t+1),  \text{ with } 0\leq d \leq t-1 \text{, and } n>k(t+1);\\[0.9em]
         \displaystyle \frac{2(n+1)}{t+1}-k & \text{ if } n \equiv t \mod(t+1)  \text{ and } n> k(t+1).\\    
        \end{cases}
        $$	
    \end{Theorem}

    \begin{proof}
   We prove the statement for $t=1$ here and for $t\geq 2$, the proof is divided into two parts, presented in Subsections \ref{subsection:lower bound} and \ref{subsection:upper bound}, where we establish the lower and upper bounds, respectively. In particular, the result follows from Theorems \ref{thm:lowerbound} and \ref{thm:upperbound}.

   For $t=1$, as mentioned after \cref{lem:dualfacets}, the ideal $I_{n,1}^{[k]}$ is the squarefree Veronese ideal of degree $k$. Therefore, $I_{n,1}^{[k]}$ is the Stanley-Reisner ideal of the $(k-2)$-dimensional skeleton of a simplex on $n$ vertices. Since the skeletons of the simplex are Cohen-Macaulay, we get $\depth(R/I_{n,1}^{[k]}) = \dim(R/I_{n,1}^{[k]}) = k-1$. Therefore, by Auslander–Buchsbaum formula, $\mathrm{pd}\left(R/I_{n,1}^{[k]}\right) = n-k+1$, as desired. 
    \end{proof}

   The following example illustrates the pattern of the projective dimension when $k$ and $t$ are fixed and the length $n$ of the path varies.

    \begin{Example} \label{Exa: proj dim formula}
    Tables \ref{Table: Proj. dim. 1} and \ref{Table: Proj. dim. 2} display the behavior of the projective dimension of $R/I_{n,4}^{[3]}$ and $R/I_{n,3}^{[4]}$ for $12\leq n\leq 29$ and $12\leq n\leq 27$, respectively.

    \begin{table}[h]
    \[
    \begin{array}{c|c|c|c|c|c|c|c|c|c|c|c|c|c|c|c|c|c|c}
    n & 12 & 13 & 14 & 15 & 16 & 17 & 18 & 19 & 20 & 21 & 22 & 23 & 24 & 25 & 26 & 27 & 28 & 29 \\[0.4em]
    \hline 
    \rule{0pt}{1.5em}
    \mathrm{pd}\left(R/I_{n,4}^{[3]}\right) & \mathbf{1} & \mathbf{2} & \mathbf{3} & 4 & 4 & 4 & 4 & \mathbf{5} & 6 & 6 & 6 & 6 & \mathbf{7} & 8 & 8 & 8 & 8 & \mathbf{9} \\
    \end{array}
    \]   
    \caption{Projective dimension of $R/I_{n,4}^{[3]}$.}
    \label{Table: Proj. dim. 1}
    \end{table}

    \begin{table}[h]
    \[
    \begin{array}{c|c|c|c|c|c|c|c|c|c|c|c|c|c|c|c|c}
    n & 12 & 13 & 14 & 15 & 16 & 17 & 18 & 19 & 20 & 21 & 22 & 23 & 24 & 25 & 26 & 27 \\
    \hline
    \rule{0pt}{1.5em}
    \mathrm{pd}\left(R/I_{n,3}^{[4]}\right) & \mathbf{1} & \mathbf{2} & \mathbf{3} & \mathbf{4} & 5 & 5 & 5 & \mathbf{6} & 7 & 7 & 7 & \mathbf{8} & 9 & 9 & 9 & \mathbf{10} \\
    \end{array}
    \]   
    \caption{Projective dimension of $R/I_{n,3}^{[4]}$.}
    \label{Table: Proj. dim. 2}
    \end{table}
    \end{Example}

   Before proceeding with the proof, we point out an interesting consequence of Theorem \ref{Theorem: Projective dimension} concerning the computation of the Leray number of $\Delta_{n,t}^{[k]}$. A simplicial complex $\Delta$ is called \emph{d-Leray} over $K$ if $\widetilde{H}_{i}(\Gamma) = 0$ for all induced subcomplexes $\Gamma \subset \Delta$ and for all $i\geq d$. The \emph{Leray number} of $\Delta$, denoted by $L(\Delta)$, is the minimum integer $d$ such that $\Delta$ is $d$-Leray over $K$. Let $K[\Delta]$ be the Stanley-Reisner ring of $\Delta$. Then, it follows from Hochster's formula that $L(\Delta) = \reg(I_\Delta)-1 = \reg(K[\Delta])$, where $\reg(I_\Delta)$ and $\reg(K[\Delta])$ denote the Castelnuovo–Mumford regularity of $I_\Delta$ and $K[\Delta]$, respectively (see \cite{KR06lerayregularity}).
    By \cite[Corollary\ 0.3]{terai99Alexanderduality}, we have $\mathrm{pd}\left(K[\Delta]\right) = \reg (K[\Delta^\vee])+1.$ 
    
    Therefore, by above discussion, we have the following corollary.
    
    \begin{Corollary}    \label{coro:leray}
        Let $n,k,t\in \NN$ such that $t\geq 1$ and  $n\geq kt$. Then, the Leray number of $\Delta_{n,t}^{[k]}$ is
    $$ 
        L(\Delta_{n,t}^{[k]})= \begin{cases}
         n-kt & \text{ if }  kt \leq n \leq k(t+1); \\[0.6em]
         \displaystyle\frac{2(n-d)}{t+1}-k & \text{ if } n \equiv d \mod(t+1),  \text{ with } 0\leq d \leq t-1 \text{, and } n>k(t+1);\\[0.9em]
         \displaystyle \frac{2(n+1)}{t+1}-k-1 & \text{ if } n \equiv t \mod(t+1)  \text{ and } n> k(t+1).\\    
        \end{cases}
        $$	
    \end{Corollary}
    \begin{proof}
          The proof follows from  the following equalities discussed above 
          \[
          \mathrm{pd}\left(R/I_{n,t}^{[k]}\right)= \mathrm{pd}\left(K[\Sigma_{n,t}^{[k]}]\right) = \reg (K[\Delta_{n,t}^{[k]}])+1 = L(\Delta_{n,t}^{[k]})+1.
          \]
    \end{proof}
    
Let us now move on to the proof of Theorem~\ref{Theorem: Projective dimension}. In the first subsection, we show that the results obtained so far allow us to provide the suitable lower bound for $\mathrm{pd}(R/I_{n,t}^{[k]})$. In general, determining the lower bound is the most challenging part from a combinatorial–algebraic perspective, since most techniques in commutative algebra provide upper bounds more readily, typically via exact sequences. We will show that the topological methods developed in the previous sections are well suited to address the problem of finding the lower bound for $\mathrm{pd}(R/I_{n,t}^{[k]})$.

\subsection{The lower bound of the projective dimension of $R/I_{n,t}^{[k]}$.}\label{subsection:lower bound}

We begin this subsection by proving some preliminary results.

\begin{Lemma}\label{prop:homology2}
    Let $n, k, t \in \NN$ such that $k> t\geq 1$ and $n = kt+l$ where $t < l \leq k$. Then
    $\widetilde{H}_{l-1}(\Delta_{n,t}^{[k]}) \neq 0$, i.e., the top homology of $\Delta_{n,t}^{[k]}$ is non-zero.    
\end{Lemma}

\begin{proof}
    We proceed by induction on $l$. Since $l>t$,  by Corollary \ref{cor:formulahomotopy}, we have the homotopy equivalence
\[
\Delta_{n,t}^{[k]} 
\simeq 
\Delta_{n-t,t}^{[k-1]} 
\vee 
\Sigma(\Delta_{n-t-1,t}^{[k-1]}) 
\vee 
\Sigma^2(\Delta_{n-t-1,t}^{[k]}).
\]

If $l= t+1$, then $n= kt+t+1$, and hence  $n-t-1 = (k-1) t+t$. 
By \cref{prop:homology}, $\Delta_{n-t-1,t}^{[k-1]} \simeq \bigvee_{{k-1 \choose t}} S^{t-1}$. Therefore, $\Sigma(\Delta_{n-t-1,t}^{[k-1]}) \simeq \bigvee_{{k-1 \choose t}} S^{t}$. Since $k-1\geq t,$ by the above homotopy decomposition of $\Delta_{n,t}^{[k]}$, it follows that $\widetilde{H}_{l-1}(\Delta_{n,t}^{[k]}) \neq 0$.

Assume that $l\geq t+2$ and that the statement holds for all smaller values of $l$. Observe that $n-t-1 = (k-1)t+(l-1)$. Then, by induction hypothesis, $\widetilde{H}_{l-2}(\Delta_{n-t-1,t}^{[k-1]})\neq 0$. 
By \cite[Theorem\ 5.16]{kozlov2007combinatorial} \footnote{There is a typo in the statement of \cite[Theorem\ 5.16]{kozlov2007combinatorial}, but the proof is correct. The correct statement is $\widetilde{H}_{l+1}(\Sigma(\Delta)) = \widetilde{H}_{l}(\Delta)$.}, $\widetilde{H}_{l-1}(\Sigma(\Delta_{n-t-1,t}^{[k-1]})) \neq 0$.
Substituting this into the above homotopy decomposition of $\Delta_{n,t}^{[k]}$ shows that  
$\widetilde{H}_{l-1}(\Delta_{n,t}^{[k]}) \neq 0$.
\end{proof}

\begin{Proposition}\label{lem:topbetti}
	Let $k,t,n\in \NN$ with $t\geq 1$ and $n\geq kt$. Then,  $\widetilde{H}_r(\Delta_{n,t}^{[k]}) \neq 0$ for 
	\[ 
	r =  \begin{cases}
		n-kt-1 & \text{ if }  kt \leq n \leq k(t+1); \\[0.6em]
		\displaystyle\frac{2n}{t+1}-k-1 & \text{ if } n>k(t+1)   \text{ and } n \equiv 0 \mod(t+1);\\[0.9em]
		\displaystyle \frac{2(n+1)}{t+1}-k-2 & \text{ if }  n>k(t+1) \text{ and } n \equiv t \mod(t+1) .\\
	\end{cases}
	\]
\end{Proposition}

\begin{proof}
If $n = kt$, then $\Delta_{n,t}^{[k]} = \{\emptyset\}$, and hence $\widetilde{H}_{-1}(\Delta_{n,t}^{[k]}) \neq 0$.  
If $kt < n \leq k(t+1)$, write $n = kt + l$ for some $1 \leq l \leq k$.  
In this case, the result follows directly from  \cref{prop:homology} and Lemma \ref{prop:homology2}.

Now assume that $n > k(t+1)$, and $n \equiv j \mod(t+1)$, where $j\in \{0,t\}$.  
Since in both cases $n-kt>t$, by Corollary \ref{cor:formulahomotopy} we have the homotopy equivalence
\[
\Delta_{n,t}^{[k]} 
\simeq 
\Delta_{n-t,t}^{[k-1]} 
\vee 
\Sigma(\Delta_{n-t-1,t}^{[k-1]}) 
\vee 
\Sigma^2(\Delta_{n-t-1,t}^{[k]}).
\]
We analyze the homology of $\Delta_{n,t}^{[k]}$ in the both congruence cases for $n$.

\textbf{Case 1:} $n>k(t+1) \text{ and } n \equiv 0 \pmod{t+1}.$  
Let $n = m(t+1)$ for some $m > k$.  
We proceed by induction on $m$.

If $m = k + 1$, then $n - t - 1 = k(t+1)$.  
By the previous case, $\widetilde{H}_{k-1}(\Delta_{n-t-1,t}^{[k]}) \neq 0$.  
By the suspension formula~\cite[Theorem\ 5.16]{kozlov2007combinatorial}, it follows that 
\[
\widetilde{H}_{k+1}(\Sigma^2(\Delta_{n-t-1,t}^{[k]})) \neq 0,
\]
and hence $\widetilde{H}_{k+1}(\Delta_{n,t}^{[k]}) \neq 0$.  
Thus, the lemma holds for $m = k+1$ because $\frac{2(k + 1)(t+1)}{t+1}-k-1 =k+1$.

Now assume $m > k + 1$ and that the statement holds for all smaller values of $m$.  
By the induction hypothesis,
\[
\widetilde{H}_{\frac{2(n - t - 1)}{t+1} - k - 1}\big(\Delta_{n-t-1,t}^{[k]}\big) \neq 0.
\]
Using the suspension formula \cite[Theorem\ 5.16]{kozlov2007combinatorial}, we get
\[
\widetilde{H}_{\left(\frac{2(n - t - 1)}{t+1} - k - 1\right) + 2}
\big(\Sigma^2(\Delta_{n-t-1,t}^{[k]})\big) \neq 0.
\]
Since $\frac{2(n - t - 1)}{t+1} - k - 1 + 2 = \frac{2n}{t+1} - k - 1$, we conclude that
\[
\widetilde{H}_{\frac{2n}{t+1} - k - 1}(\Delta_{n,t}^{[k]}) \neq 0.
\]
This establishes the result in this case.

\textbf{Case 2:} $n>k(t+1) \text{ and } n \equiv t \pmod{t+1}.$  
Write $n = m(t+1) + t$ for some $m \geq k$. We again proceed by induction on $m$.  

If $m = k$, then $n - t - 1 = kt + k - 1$. 
By the earlier case, $\widetilde{H}_{k-2}(\Delta_{n-t-1,t}^{[k]}) \neq 0$.  
Therefore,
\[
\widetilde{H}_{k}(\Sigma^2(\Delta_{n-t-1,t}^{[k]})) \neq 0,
\]
and consequently $\widetilde{H}_{k}(\Delta_{n,t}^{[k]}) \neq 0$.  
Thus, the lemma holds for $m = k$ because $\frac{2(k(t+1)+t+1)}{t+1}-k-2 = k$.

Now suppose $m > k$ and that the statement holds for all smaller values of $m$.  
By the induction hypothesis,
\[
\widetilde{H}_{\frac{2(n - t - 1 + 1)}{t+1} - k - 2}
\big(\Delta_{n-t-1,t}^{[k]}\big) \neq 0.
\]
Using the suspension formula, we obtain
\[
\widetilde{H}_{\left(\frac{2(n - t - 1 + 1)}{t+1} - k - 2\right) + 2}
\big(\Sigma^2(\Delta_{n-t-1,t}^{[k]})\big) \neq 0.
\]
Since $(\frac{2(n - t - 1 + 1)}{t+1} - k - 2) + 2 = \frac{2(n + 1)}{t+1} - k - 2$, we conclude that
\[
\widetilde{H}_{\frac{2(n + 1)}{t+1} - k - 2}(\Delta_{n,t}^{[k]}) \neq 0.
\]
This completes the proof.
\end{proof}
\begin{Theorem}\label{thm:lowerbound}
	Under the Notation \ref{Notation}, for any $1\leq k \leq \nu(\Gamma_{n,t})$, we have
	
    \[
\mathrm{pd}\left(R/I_{n,t}^{[k]}\right) \geq 
\begin{cases}
n - kt + 1, 
& \text{if } kt \leq n \leq k(t+1); \\[0.6em]
\displaystyle \frac{2(n - d)}{t+1} - k + 1, 
& \text{if } n \equiv d \pmod{t+1},\ 
0 \leq d \leq t-1,\ \text{and } n > k(t+1); \\[0.9em]
\displaystyle \frac{2(n + 1)}{t+1} - k, 
& \text{if } n \equiv t \pmod{t+1},\ \text{and } n > k(t+1).
\end{cases}
\]
\end{Theorem}

\begin{proof}
It suffices to show that there exists a nonzero Betti number $\beta_{r,m}$ for some $m \geq 0$, where $r$ is given by the corresponding expression on the right-hand side above.
Applying Hochster’s formula (see \cite[Corollary\ 8.1.4]{HHmonomialideals}) for $R/I_{n,t}^{[k]}$, we get
    \[
    \beta_{i,n}(R/I_{n,t}^{[k]}) = \dim_K \widetilde{H}_{i-2}(\lk_{\Delta_{n,t}^{[k]}}(\emptyset); K)
    = \dim_K \widetilde{H}_{i-2}(\Delta_{n,t}^{[k]}; K).
    \]

\textbf{Case 1:} When $kt \leq n \leq k(t+1)$ or  $n > k(t+1)$ and $n \equiv 0 \pmod{t+1}$ or $n > k(t+1)$ and $n \equiv t \pmod{t+1}$.
Then, the result follows from Proposition \ref{lem:topbetti}.

\medskip
\noindent\textbf{Case 2.}
Suppose that $n>k(t+1)$ and
$n\equiv d\pmod{t+1}$ for some $1\leq d\leq t-1$.
Set $W=[n-d]$ and $S=K[x_i:i\in W]$. We first observe that $I_{n,t}^{[k]}\cap S=I_{n-d,t}^{[k]}.$
It follows that
\[
    \Sigma_{n-d,t}^{[k]}
    =
    (\Sigma_{n,t}^{[k]})_W.
\]

For a squarefree multidegree $F$, let $\beta_{i,F}(M)$ denote the
corresponding multigraded Betti number of $M$ over its ambient
polynomial ring. By Hochster's formula and the equality above, we
have
\[
   \beta_{i,W}
    \left(S/I_{n-d,t}^{[k]}\right)
    =
    \dim_K
    \widetilde H_{|W|-i-1}
    \left(\Sigma_{n-d,t}^{[k]};K\right)\\
    =
    \dim_K
    \widetilde H_{|W|-i-1}
    \left(
        \left.\Sigma_{n,t}^{[k]}\right|_W;K
    \right)\\
    =
    \beta_{i,W}
    \left(R/I_{n,t}^{[k]}\right).
\]
Set $r=\frac{2(n-d)}{t+1}-k+1.$ By Case~1, $\beta_{r,n-d}
    \left(S/I_{n-d,t}^{[k]}\right)\neq0.$ Since $W$ is the unique subset of $W$ having cardinality $n-d$, we have $\beta_{r,n-d}
    \left(S/I_{n-d,t}^{[k]}\right)
    =
    \beta_{r,W}
    \left(S/I_{n-d,t}^{[k]}\right).$ Therefore, $\beta_{r,W}
    \left(R/I_{n,t}^{[k]}\right)\neq0.$
Since a multigraded Betti number is a summand of the corresponding
graded Betti number, it follows that $\beta_{r,n-d}
    \left(R/I_{n,t}^{[k]}\right)\neq0,$ that is,
\[
    \operatorname{pd}
    \left(R/I_{n,t}^{[k]}\right)
    \geq r=
    \frac{2(n-d)}{t+1}-k+1,
\]
as desired.

\end{proof}

\subsection{The upper bound of the projective dimension of $R/I_{n,t}^{[k]}$}\label{subsection:upper bound}

This subsection is devoted to establishing an upper bound for the projective dimension of $R/I_{n,t}^{[k]}$. We begin by introducing the function $\chi_t(n,k)$, which anticipates the formulas that will later be proved. For every integer $n$ and all positive integers $k$ and $t$, we define
$$ 
\chi_t(n,k) =  \begin{cases}
        0 & \text{ if }   n <kt; \\[0.6em]
	n-kt+1 & \text{ if }  kt \leq n \leq k(t+1); \\[0.6em]
	\displaystyle\frac{2(n-d)}{t+1}-k+1 & \text{ if } n \equiv d \mod(t+1),  \text{ with } 0\leq d \leq t-1 \text{ and } n>k(t+1);\\[0.9em]
	\displaystyle \frac{2(n+1)}{t+1}-k & \text{ if } n \equiv t \mod(t+1),  \text{ with } n>k(t+1).\\
\end{cases}
$$	

To establish the upper bound in Theorem~\ref{Theorem: Projective dimension}, our approach relies on certain monotonicity and comparison properties of $\chi_t(n,k)$.  
The following lemma provides one such comparison, relating the values of $\chi_t$ when both $n$ and $k$ are simultaneously reduced.

\begin{Lemma}\label{lem:pdrelation}
With the notation above, we have $\chi_t(n-it,k+1 -i) \leq \chi_t(n,k+1)$ for all $1\leq i \leq k$.
\end{Lemma}     

\begin{proof}
We begin by claiming that for all  $0\leq i\leq k-1$,
\[
\chi_t(n-(i+1)t,k+1 -(i+1))\leq \chi_t(n-it,k+1 -i).
\]
Once this claim is established, it immediately follows that
\[
\chi_t(n-kt,1)\leq \chi_t(n-(k-1)t,2)\leq \chi_t(n-(k-2)t,3)\leq \ldots\leq \chi_t(n-t,k) \leq \chi_t(n,k+1),
\]
which yields the desired inequality.

To prove the claim, fix an integer $0 \leq i \leq k-1$, and set $m=k+1-i$. We consider four cases according to the value of $n$. 

\medskip
\noindent\textbf{Case 1.}
Suppose that $n-it<mt$. Then
$\chi_t(n-it,m)=0$. Moreover,
\[
    n-(i+1)t
    <
    mt-t
    =
    (m-1)t,
\]
and therefore
$\chi_t(n-(i+1)t,m-1)=0$. Hence
\[
    \chi_t(n-(i+1)t,m-1)
    =
    \chi_t(n-it,m)
    =
    0.
\]

\noindent\textbf{Case 2.} Suppose that $mt  \leq n-it \leq m(t+1)$. The pair $(n-it,m)$ belong to the second case in the definition of $\chi_t$, and we have $\chi_t(n-it,m) = (n-it) - mt + 1.$
Subtracting $t$ from each term of $mt  \leq n-it \leq m(t+1)$ gives
        \[	
	(m-1)t \leq n-(i+1)t \leq (m-1)(t+1) +1. 
        \]
If $n-(i+1)t \leq (m-1)(t+1)$, then the pair \((n-(i+1)t,m-1)\) also belong to the second case in the definition of $\chi_t$, and we obtain
\begin{align*}
\chi_t(n-(i+1)t,m-1) = (n-(i+1)t) - (m-1)t + 1= n - it - mt + 1,
\end{align*}
which coincides with  $\chi_t(n-it,m)$, as claimed.

Now, let $n-(i+1)t = (m-1)(t+1)+1$. Then $n-(i+1)t \equiv 1 \pmod{t+1}$. 
Hence, by the third case in the definition of $\chi_t$, we have
\begin{align*}
\chi_t(n-(i+1)t,m-1) 
&= \frac{2\big(n-(i+1)t-1\big)}{t+1} - (m-1) + 1. 
\end{align*}    
Then, substituting $n-(i+1)t-1 = (m-1)(t+1)$ in above equality gives $\chi_t(n-(i+1)t,m-1)  = m.$
It is straightforward to check that  $m = n-it-mt$. Therefore,
\[
\chi_t(n-(i+1)t,m-1) 
= \chi_t(n-it,m) - 1,
\]
as claimed.

\noindent\textbf{Case 3.} Suppose $n-it > m(t+1)$ and $n-it \equiv d \pmod{t+1}$ with $0 \leq d \leq t-1$. Under these conditions, $\chi_t(n - it, m)$ falls under the third case of its definition. Note that $n-(i+1)t \equiv d+1 \pmod{t+1}$, and $n-(i+1)t>(m-1)(t+1)$.  

\noindent If $1 \leq d+1 \leq t-1$, then 
$$
\chi_t(n-(i+1)t,m-1) = \frac{2\big(n-(i+1)t-(d+1)\big)}{t+1} - (m-1) + 1 = \chi_t(n-it,m) - 1,
$$
and if $d+1 = t$, then 
$$\chi_t\big(n-(i+1)t,m-1\big)   =
	\frac{2\big(n-(i+1)t+1 \big)}{t+1} -(m-1) 	 = \chi_t(n-it,m)
$$
as claimed. 

\noindent\textbf{Case 4.} Suppose $n-it > m(t+1)$ and $n-it \equiv t \pmod{t+1}$. Then $n-(i+1)t \equiv 0 \pmod{t+1}.$ Under these conditions, $\chi_t(n - it, m)$ falls under the fourth case and $\chi_t(n-(i+1)t,m-1)$ falls under the third case of its definition. This gives
$$
\chi_t(n-(i+1)t,m-1) = \frac{2(n-(i+1)t)}{t+1} - (m-1) + 1=\chi_t(n - it, m).
$$
as claimed.
\end{proof}

The next lemma is to understand the behavior of $\chi_t(n,k)$ as a function of $n$ for fixed $t$ and $k$. We show that it is an non-decreasing function of $n$. We also give quantitative control on $\chi_t(n,k)$ when $n$ decreases by $t$ or $t+1$.

\begin{Lemma}\label{lem:pdincreasing}
	For $t\geq 1$ and $k\geq 1$, $\chi_t(n,k)$ is a non-decreasing function of $n$. Moreover, 
\[
\chi_t(n-t,k)<\chi_t(n,k) \qand \chi_t(n-(t+1),k)\leq \chi_t(n,k)-2 \qforall n> kt.
\]
\end{Lemma}

\begin{proof}
For every integer $m$, set
$\delta_k(m)=\chi_t(m+1,k)-\chi_t(m,k)$.
A direct computation from the definition of $\chi_t$ gives
\[
\delta_k(m)=
\begin{cases}
0,
    & m<kt-1,\\
1,
    & kt-1\leq m<k(t+1),\\
1,
    & m\geq k(t+1)
      \text{ and }
      m\bmod(t+1)\in\{t-1,t\},\\
0,
    & m\geq k(t+1)
      \text{ and }
      m\bmod(t+1)\notin\{t-1,t\}.
\end{cases}
\]
In particular, $\delta_k(m)\geq0$ for every integer $m$, and hence
$\chi_t(n,k)$ is non-decreasing.

We next prove the first strict inequality. If $n-t<kt$, then
$\chi_t(n-t,k)=0$, whereas $\chi_t(n,k)>0$ because $n>kt$.
Suppose now that $n-t\geq kt$. Then
\[
    \chi_t(n,k)-\chi_t(n-t,k)
    =
    \sum_{m=n-t}^{n-1}\delta_k(m).
\]
The set $\{n-t,n-t+1,\dots,n-1\}$ consists of $t$ consecutive integers.
Moreover, observe that this set contains at least one integer $m$ for which $\delta_k(m)=1$. Therefore, $\chi_t(n-t,k)<\chi_t(n,k)$.

Finally, suppose first that $n-(t+1)<kt$. Then $\chi_t(n-(t+1),k)=0$. Since $n>kt$, we have
$\chi_t(n,k)\geq2$, and hence \[
    \chi_t(n-(t+1),k)
    \leq
    \chi_t(n,k)-2.
\]
Suppose now that $n-(t+1)\geq kt$. Then
\[
    \chi_t(n,k)-\chi_t(n-(t+1),k)
    =
    \sum_{m=n-(t+1)}^{n-1}\delta_k(m).
\]
The set $\{n-(t+1),n-t,\dots,n-1\}$ consists of $t+1$ consecutive integers and it contains at least two integers $m$ for which $\delta_k(m)=1$. Consequently,
\[
    \chi_t(n-(t+1),k)
    \leq
    \chi_t(n,k)-2.
\]
\end{proof}

We next compare values of $\chi_t$ when $n$ is shifted by $(t+1)$ but $k$ is reduced by $1$. The following inequality is crucial in the proof of our main theorem.

\begin{Lemma}\label{lem:sub_t-1}
$\chi_t(n-(t+1),k) \leq \chi_t(n,k+1)-1$ for $t \geq 2,$ $k \geq 1$ and $n\geq (k+1)t$. 
\end{Lemma}

\begin{proof}
Since $n \geq (k + 1)t$, it suffices to consider the following cases. 

\noindent\textbf{Case 1.} Suppose $(k + 1)t \leq n \leq (k + 1)(t + 1)$. Then $kt - 1 \leq n - t - 1 \leq k(t + 1)$.

If $n-t-1 = kt-1$, then $ \chi_t(n-t-1,k) = 0 = \chi_t(n,k+1)-1.$ On the other hand, if $kt \leq n-t-1 \leq k(t+1)$, then by the second piece of the definition,
\begin{align*}
\chi_t(n-t-1,k) 
&= (n-t-1) - kt + 1 = n-(k+1)t+1-1 = \chi_t(n,k+1)-1.
\end{align*}
	
\noindent\textbf{Case 2.} Suppose $n > (k+1)(t+1)$ and $n \equiv d \pmod{t+1}$ with $0 \leq d \leq t$.  Then $n-(t+1) \equiv d \pmod{t+1}$. By the definition of $\chi_t$ (the third case if $d \neq t$ and the fourth case if $d = t$), we have
\[
\chi_t(n - t - 1, k) = \chi_t(n, k + 1) - 1,
\]
as required. 
\end{proof}

\begin{Lemma}\label{lem:pdequation}
Let $t\geq 2$, $k\geq 1$, $n>t(k+1)$ and
$2\leq i\leq n-t(k+1)+1$. Then
\[
    \chi_t(n-t-i+1,k)+\chi_t(i-2,1)
    \leq \chi_t(n-t-1,k).
\]
\end{Lemma}

\begin{proof}

Set $b=n-t-i+1$ and $a=i-2$.
Since $i\leq n-t(k+1)+1$, we have
$b=n-t-i+1\geq n-t-\bigl(n-t(k+1)+1\bigr)+1=kt$.
Moreover, $a+b=(i-2)+(n-t-i+1)=n-t-1$.

Write $a=q(t+1)+r$, where $q\geq0$ and $0\leq r\leq t$.
By the definition of $\chi_t(\,\cdot\,,1)$, we have
\[
    \chi_t(a,1)=
    \begin{cases}
        2q,   & 0\leq r\leq t-1,\\
        2q+1, & r=t.
    \end{cases}
\]

Since $b\geq kt$, Lemma~\ref{lem:pdincreasing} implies
$\chi_t(m+t+1,k)\geq\chi_t(m,k)+2$ for every $m\geq kt$.
Consequently, applying this inequality successively $q$ times gives
\begin{equation}\label{eq:q-blocks}
    \chi_t\bigl(b+q(t+1),k\bigr)
    \geq \chi_t(b,k)+2q.
\end{equation}

Suppose first that $0\leq r\leq t-1$. Since
$\chi_t(\,\cdot\,,k)$ is non-decreasing, it follows from
\eqref{eq:q-blocks} that
\[
    \chi_t(a+b,k)=\chi_t\bigl(b+q(t+1)+r,k\bigr)
    \geq\chi_t\bigl(b+q(t+1),k\bigr)
    \geq\chi_t(b,k)+2q
    =\chi_t(b,k)+\chi_t(a,1).
\]

Suppose now that $r=t$. Again by
Lemma~\ref{lem:pdincreasing}, we have
$\chi_t(m+t,k)>\chi_t(m,k)$ for every $m\geq kt$.
Since $\chi_t$ takes integer values, this implies
$\chi_t(m+t,k)\geq\chi_t(m,k)+1$.
Using this inequality with $m=b+q(t+1)$ and then
\eqref{eq:q-blocks}, we obtain
\[    \chi_t(a+b,k)
    =\chi_t\bigl(b+q(t+1)+t,k\bigr)
    \geq\chi_t\bigl(b+q(t+1),k\bigr)+1
    \geq\chi_t(b,k)+2q+1
    =\chi_t(b,k)+\chi_t(a,1).
\]

Since $a+b=n-t-1$, $a=i-2$, and $b=n-t-i+1$, the
inequality obtained above becomes $\chi_t(n-t-i+1,k)+\chi_t(i-2,1)
    \leq\chi_t(n-t-1,k),$ as desired.
\end{proof}






We are now ready to prove the upper bound of the projective dimension of $R/I_{n,t}^{[k]}$.

\begin{Theorem}\label{thm:upperbound}
	Under the Notation \ref{Notation}, for any $1\leq k \leq \nu(\Gamma_{n,t})$ and $t\geq 2$, we have
	$$
\mathrm{pd}\!\left(R/I_{n,t}^{[k]}\right)\leq \chi_t(n, k) =
\begin{cases}
n-kt+1 & \text{if } kt \le n \le k(t+1); \\[0.8em]
\displaystyle \frac{2(n-d)}{t+1}-k+1 &
\text{if }
\text{
$\begin{array}{l}
n \equiv d \pmod{t+1}, \text{ with } 0 \le d \le t-1,\\
\text{and } n > k(t+1);
\end{array}$
} \\[1.3em]
\displaystyle \frac{2(n+1)}{t+1}-k &
\text{if }
 n \equiv t \pmod{t+1} \text{ and } n > k(t+1).
\end{cases}
$$

\end{Theorem}

\begin{proof}
We proceed by induction on $n$. 

If $n=t$, then $\nu(\Gamma_{n,t})=1$ and
$I_{t,t}=(x_1x_2\cdots x_t)$. Therefore,
$
    \operatorname{pd}(R/I_{t,t})
    =1
    =
    \chi_t(t,1).
$

Now let $n>t$, and assume that the assertion holds for every path
on $m$ vertices with $t\leq m<n$ and every
$1\leq\ell\leq\nu(\Gamma_{m,t})$. Fix
$1\leq\ell\leq\nu(\Gamma_{n,t})$.

If $\ell=1$, the result follows from
\cite[Theorem~4.1]{HT2010pathideals}. We may therefore assume that
$\ell=k+1$ for some $k\geq1$. We prove that
\[
    \operatorname{pd}
    \left(R/I_{n,t}^{[k+1]}\right)
    \leq
    \chi_t(n,k+1).
\]

Set $a=n-t(k+1)$. Since
$k+1\leq\nu(\Gamma_{n,t})$, we have $a\geq0$.
We first consider the case $a=0$. In this case
$n=t(k+1)$ and
$
    I_{n,t}^{[k+1]}=(x_1x_2\cdots x_n).
$, consequently,
$
    \operatorname{pd}
    \left(R/I_{n,t}^{[k+1]}\right)
    =1
    =
    \chi_t(n,k+1).
$
Thus, in the remainder of the proof, we may assume that $a\geq1$,
or equivalently, that $n>t(k+1)$.

For integers $1\leq p\leq q\leq n$, let $P_{[p,q]}$ denote the path
graph on the vertex set $[p,q]:=\{p,p+1,\ldots,q\}$, and let
$I_{[p,q]}$ denote its $t$-path ideal. For
$1\leq i\leq n-t+1$, set
$f_i=x_i x_{i+1}\cdots x_{i+t-1}$. Thus,
\[
    I_{[p,q]}=(f_p,f_{p+1},\ldots,f_{q-t+1})
\]
whenever $q-p+1\geq t$, and we set $I_{[p,q]}=0$ when
$q-p+1<t$.

Consider the exact sequence
\[
0\longrightarrow
\frac{R}{I_{[1,n]}^{[k+1]}:f_1}(-t)
\xrightarrow{\ \cdot f_1\ }
\frac{R}{I_{[1,n]}^{[k+1]}}
\longrightarrow
\frac{R}{I_{[1,n]}^{[k+1]}+(f_1)}
\longrightarrow0.
\]
For $1\leq i\leq a+1$, the standard colon computations give
\[
    I_{[i,n]}^{[k+1]}:f_i
    =
    I_{[i+t,n]}^{[k]}
\]
and, for $2\leq i\leq a+1$,
\[
    I_{[1,t+i-2]}:f_i
    =(f_1,\ldots,f_{i-1}):f_i
    =I_{[1,i-2]}+(x_{i-1}).
\]
Moreover,
\[
\left(
I_{[i,n]}^{[k+1]}+I_{[1,t+i-2]}
\right):f_i=
I_{[i+t,n]}^{[k]}+
\left(I_{[1,t+i-2]}:f_i\right),
\]
whereas
\[
I_{[i,n]}^{[k+1]}+I_{[1,t+i-2]}+(f_i)=
I_{[i+1,n]}^{[k+1]}+I_{[1,t+i-1]}.
\]
Since $I_{[a+2,n]}^{[k+1]}=0$, successive applications of the
preceding exact sequence,
yield
\[
\operatorname{pd}
\left(R/I_{n,t}^{[k+1]}\right)
\leq
\max
\left\{
\operatorname{pd}
\left(R/I_{[1+t,n]}^{[k]}\right),
\operatorname{pd}
\left(R/I_{[1,t+a]}\right),
\alpha
\right\},
\]
where
\[
\alpha=
\max_{2\leq i\leq a+1}
\left\{
\operatorname{pd}
\left(
\frac{R}{
I_{[i+t,n]}^{[k]}+
(I_{[1,t+i-2]}:f_i)}
\right)
\right\}.
\]

For $2\leq i\leq a+1$, the ideals
$I_{[i+t,n]}^{[k]}$ and $I_{[1,t+i-2]}:f_i$ belong to $K[x_{i+t},\ldots,x_n]$ and $K[x_1,\ldots,x_{i-1}]$, respectively.
Thus, tensoring their minimal free resolutions gives
\[
\operatorname{pd}
\left(
\frac{R}{
I_{[i+t,n]}^{[k]}+
(I_{[1,t+i-2]}:f_i)}
\right)
=
\operatorname{pd}
\left(R/I_{[i+t,n]}^{[k]}\right)\\
\quad+
\operatorname{pd}
\left(
R/(I_{[1,t+i-2]}:f_i)
\right).
\]
Furthermore,
\[
\operatorname{pd}
\left(
R/(I_{[1,t+i-2]}:f_i)
\right)
=
\operatorname{pd}
\left(R/I_{[1,i-2]}\right)+1.
\]

We now estimate the terms in the maximum. The path
$P_{[1+t,n]}$ has $n-t<n$ vertices. Hence, by the induction
hypothesis,
\[
\operatorname{pd}
\left(R/I_{[1+t,n]}^{[k]}\right)
\leq
\chi_t(n-t,k)
\leq
\chi_t(n,k+1),
\]
where the last inequality follows from
Lemma~\ref{lem:pdrelation}.

The path $P_{[1,t+a]}$ has $t+a=n-kt<n$ vertices. Therefore,
again by the induction hypothesis,
\[
\operatorname{pd}
\left(R/I_{[1,t+a]}\right)
\leq
\chi_t(n-kt,1)
\leq
\chi_t(n,k+1),
\]
where the last inequality follows from
Lemma~\ref{lem:pdrelation}.

Finally, let $2\leq i\leq a+1$. The path
$P_{[i+t,n]}$ has $n-t-i+1<n$ vertices. Moreover,
$i\leq a+1$ implies $n-t-i+1\geq kt$. Thus, the induction
hypothesis gives
\[
    \operatorname{pd}
    \left(R/I_{[i+t,n]}^{[k]}\right)
    \leq
    \chi_t(n-t-i+1,k).
\]
Similarly,
\[
    \operatorname{pd}
    \left(R/I_{[1,i-2]}\right)
    \leq
    \chi_t(i-2,1),
\]
where, if $i-2<t$, both sides are zero. Consequently,
\begin{align*}
\operatorname{pd}
\left(
\frac{R}{
I_{[i+t,n]}^{[k]}+
(I_{[1,t+i-2]}:f_i)}
\right)
&\leq
\chi_t(n-t-i+1,k)+\chi_t(i-2,1)+1\\
&\leq
\chi_t(n-t-1,k)+1
\qquad\text{by Lemma~\ref{lem:pdequation}}\\
&\leq
\chi_t(n,k+1)
\qquad\text{by Lemma~\ref{lem:sub_t-1}}.
\end{align*}

Since every term in the maximum is bounded above by
$\chi_t(n,k+1)$, we conclude that
\[
    \operatorname{pd}
    \left(R/I_{n,t}^{[k+1]}\right)
    \leq
    \chi_t(n,k+1).
\]
This completes the proof.
\end{proof}

\subsection{Krull dimension of $R/I_{n,t}^{[k]}$.}\label{subsection: Krull}
We conclude the paper with this short subsection, in which we give a combinatorial formula for the Krull dimension of $R/I_{n,t}^{[k]}$.

\begin{Theorem}\label{Theorem: Krull dimension}
        Under the Notation \ref{Notation}, for any $1\leq k \leq \nu(\Gamma_{n,t})$, we have
        $$ 
        \mathrm{dim}\left(R/I_{n,t}^{[k]}\right)=n-\left\lfloor \frac{n}{t}\right\rfloor+k-1.
        $$	
    \end{Theorem}

\begin{proof}


To prove the assertion, it is enough to show that 
$$ 
 	\mathrm{height}\left(I_{n,t}^{[k]}\right)=\nu(\Gamma_{n,t})-k+1.
$$	

Set $J=I_{n,t}^{[k]}$ and
$\nu=\nu(\Gamma_{n,t})=\left\lfloor n/t\right\rfloor$. For $1\leq j\leq\nu$, let
$B_j=\{(j-1)t+1,\ldots,jt\}$.
These are pairwise disjoint $t$-paths in $P_n$.

Set $C_0=\{t,2t,\ldots,(\nu-k+1)t\}$. We claim that $C_0$
intersects the support of every generator of $J$. Indeed, let
$u\in G(J)$. Then $\operatorname{supp}(u)$ is the union of $k$
pairwise disjoint intervals of length $t$. Each such interval contains
exactly one multiple of $t$, and the multiples obtained from distinct
intervals are distinct. Since only $k-1$ elements of
$\{t,2t,\ldots,\nu t\}$ lie outside $C_0$, at least one of these
$k$ multiples belongs to $C_0$. Consequently,
$\operatorname{supp}(u)\cap C_0\neq\emptyset$, and hence
$\height(J)\leq |C_0|=\nu-k+1$.

Conversely, let $C\subseteq[n]$ intersect the support of every
generator of $J$. If $C$ intersected at most $\nu-k$ of the sets $B_1,\ldots,B_\nu$, then at least $k$ of these sets would be
disjoint from $C$. The union of any such $k$ sets supports a
generator of $J$, contradicting the defining property of $C$.
Therefore, $C$ intersects at least $\nu-k+1$ sets of $B_1,\ldots,B_\nu$, which are disjoints. Then $|C|\geq\nu-k+1$. Thus $\height(J)\geq\nu-k+1$. In conclusion, $\height\left(I_{n,t}^{[k]}\right)=\nu-k+1$.
\end{proof}

	
    

    \begin{Example}

    Tables \ref{Table: Krull dim. 1} and \ref{Table: Krull dim. 2} illustrate the pattern observed for the Krull dimension of $R/I_{n,2}^{[2]}$ and $R/I_{n,5}^{[3]}$ when $4\leq n\leq 19$ and $15\leq n\leq 29$, respectively.

    \begin{table}[h]
    \[
    \begin{array}{c|c|c|c|c|c|c|c|c|c|c|c|c|c|c|c|c}
    n & 4 & 5 & 6 & 7 & 8 & 9 & 10 & 11 & 12 & 13 & 14 & 15 & 16 & 17 & 18 & 19  \\
    \hline
    \rule{0pt}{1.5em}
    \mathrm{dim}\left(R/I_{n,2}^{[2]}\right) & \mathbf{3} & \mathbf{4} & 4 & 5 & \mathbf{5} & \mathbf{6} & 6 & 7 & \mathbf{7} & \mathbf{8} & 8 & 9 & \mathbf{9} & \mathbf{10} & 10 & 11 \\
    \end{array}
    \]  
    \caption{Krull dimension of $R/I_{n,2}^{[2]}$.}
    \label{Table: Krull dim. 1}
    \end{table}

    \begin{table}[h]
    \[
    \begin{array}{c|c|c|c|c|c|c|c|c|c|c|c|c|c|c|c}
    n &  15 & 16 & 17 & 18 & 19 & 20 & 21 & 22 & 23 & 24 & 25 & 26 & 27 & 28 & 29 \\[0.4em]
    \hline 
    \rule{0pt}{1.5em}
    \mathrm{dim}\left(R/I_{n,5}^{[3]}\right) &  \mathbf{14} & \mathbf{15} & \mathbf{16} & \mathbf{17} & \mathbf{18} & 18 & 19 & 20 & 21 & 22 & \mathbf{22} & \mathbf{23} & \mathbf{24} & \mathbf{25} & \mathbf{26} \\
    \end{array}
    \]   
    \caption{Krull dimension of $R/I_{n,5}^{[3]}$.}
    \label{Table: Krull dim. 2}
    \end{table}
    
    \end{Example}

     


\begin{thebibliography}{EHHM24}

\bibitem[AF15]{AF}
Ali Alilooee and Sara Faridi.
\newblock On the resolution of path ideals of cycles.
\newblock {\em Communications in Algebra}, 43(12):5413--5433, 2015.

\bibitem[AF18]{AF2}
A.~Alilooee and S.~Faridi.
\newblock Graded betti numbers of path ideals of cycles and lines.
\newblock {\em Journal of Algebra and Its Applications}, 17(01):1850011, 2018.

\bibitem[BCV25]{BCV}
Silviu Bălănescu, Mircea Cimpoeaş, and Thanh Vu.
\newblock Betti numbers of powers of path ideals of cycles.
\newblock {\em Journal of Algebraic Combinatorics}, 61(4):44, 2025.

\bibitem[Ber89]{B}
C.~Berge.
\newblock {\em Hypergraphs: Combinatorics of Finite Sets}, volume~45 of {\em Mathematical Library}.
\newblock North-Holland, 1989.

\bibitem[BH98]{BH}
Winfried Bruns and H.~Jürgen Herzog.
\newblock {\em Cohen-Macaulay Rings}.
\newblock Cambridge Studies in Advanced Mathematics. Cambridge University Press, 2 edition, 1998.

\bibitem[BHO11]{T}
Rachelle~R. Bouchat, Huy~Tài Hà, and Augustine OʼKeefe.
\newblock Path ideals of rooted trees and their graded betti numbers.
\newblock {\em Journal of Combinatorial Theory, Series A}, 118(8):2411--2425, 2011.

\bibitem[BHZN18]{BHZ18powersedgeideals}
Mina Bigdeli, J\"urgen Herzog, and Rashid Zaare-Nahandi.
\newblock On the index of powers of edge ideals.
\newblock {\em Comm. Algebra}, 46(3):1080--1095, 2018.

\bibitem[BJMV23]{Bayer23matching}
Margaret Bayer, Marija Jeli\'c{}~Milutinovi\'c, and Julianne Vega.
\newblock General polygonal line tilings and their matching complexes.
\newblock {\em Discrete Math.}, 346(7):Paper No. 113428, 12, 2023.

\bibitem[BW96]{BW}
Anders Björner and Michelle~L. Wachs.
\newblock Shellable nonpure complexes and posets.
\newblock {\em Transactions of the American Mathematical Society}, 348(4):1299--1327, 1996.

\bibitem[BW97]{BW97part2}
Anders Bj{\"o}rner and Michelle~L. Wachs.
\newblock Shellable nonpure complexes and posets. {II}.
\newblock {\em Trans. Amer. Math. Soc.}, 349(10):3945--3975, 1997.

\bibitem[BWW09]{BWV09SCM}
Anders Bj\"orner, Michelle Wachs, and Volkmar Welker.
\newblock On sequentially {C}ohen-{M}acaulay complexes and posets.
\newblock {\em Israel J. Math.}, 169:295--316, 2009.

\bibitem[CD99]{CN}
Aldo Conca and Emanuela {De Negri}.
\newblock M-sequences, graph ideals, and ladder ideals of linear type.
\newblock {\em Journal of Algebra}, 211(2):599--624, 1999.

\bibitem[CFL25]{CFL}
Marilena Crupi, Antonino Ficarra, and Ernesto Lax.
\newblock {Matchings, squarefree powers, and Betti splittings}.
\newblock {\em Illinois Journal of Mathematics}, 69(2):353 -- 372, 2025.

\bibitem[DRS24]{Kamalesh2}
Kanoy~Kumar Das, Amit Roy, and Kamalesh Saha.
\newblock Square-free powers of {C}ohen-{M}acaulay forests, cycles, and whiskered cycles.
\newblock {\em arXiv:2409.06021}, 2024.

\bibitem[DRS26]{Kamalesh1}
Kanoy~Kumar Das, Amit Roy, and Kamalesh Saha.
\newblock Square-free powers of {C}ohen-{M}acaulay simplicial forests.
\newblock {\em Proc. Amer. Math. Soc.}, 154(2):549--565, 2026.

\bibitem[EF26]{EF}
Nursel Erey and Antonino Ficarra.
\newblock Matching powers of monomial ideals and edge ideals of weighted oriented graphs.
\newblock {\em Journal of Algebra and Its Applications}, 25(11):2650118, 2026.

\bibitem[EH21]{EH}
N.~Erey and T.~Hibi.
\newblock Squarefree powers of edge ideals of forests.
\newblock {\em Electronic Journal of Combinatorics}, 28(2):P2.32, 2021.

\bibitem[EHHM22]{EHHM}
N.~Erey, J.~Herzog, T.~Hibi, and S.~S. Madani.
\newblock Matchings and squarefree powers of edge ideals.
\newblock {\em Journal of Combinatorial Theory, Series A}, 188:105585, 2022.

\bibitem[EHHM24]{EHHM2}
N.~Erey, J.~Herzog, T.~Hibi, and S.~S. Madani.
\newblock The normalized depth function of squarefree powers.
\newblock {\em Collectanea Mathematica}, 75:409--423, 2024.

\bibitem[ER98]{ER98linearresolution}
John~A. Eagon and Victor Reiner.
\newblock Resolutions of {S}tanley-{R}eisner rings and {A}lexander duality.
\newblock {\em J. Pure Appl. Algebra}, 130(3):265--275, 1998.

\bibitem[Fak24]{S}
S.~A.~Seyed Fakhari.
\newblock On the castelnuovo–mumford regularity of squarefree powers of edge ideals.
\newblock {\em Journal of Pure and Applied Algebra}, 228(3):107488, 2024.

\bibitem[FHH23]{FHH}
A.~Ficarra, J.~Herzog, and T.~Hibi.
\newblock Behaviour of the normalized depth function.
\newblock {\em Electronic Journal of Combinatorics}, 30(2):P2.31, 2023.

\bibitem[GS]{M2}
Daniel~R. Grayson and Michael~E. Stillman.
\newblock Macaulay2, a software system for research in algebraic geometry.
\newblock Available at \url{http://www.math.uiuc.edu/Macaulay2/}.

\bibitem[Hat02]{hatcher}
Allen Hatcher.
\newblock {\em Algebraic topology}.
\newblock Cambridge University Press, Cambridge, 2002.

\bibitem[HH11]{HHmonomialideals}
J\"urgen Herzog and Takayuki Hibi.
\newblock {\em Monomial ideals}, volume 260 of {\em Graduate Texts in Mathematics}.
\newblock Springer-Verlag London, Ltd., London, 2011.

\bibitem[HHTZ08]{HHTZ}
J.~Herzog, T.~Hibi, N.~V. Trung, and X.~Zheng.
\newblock Standard graded vertex cover algebras, cycles and leaves.
\newblock {\em Transactions of the American Mathematical Society}, 360(12):6231--6249, 2008.

\bibitem[Hoc77]{hochsterformula}
Melvin Hochster.
\newblock Cohen-{M}acaulay rings, combinatorics, and simplicial complexes.
\newblock In {\em Ring theory, {II} ({P}roc. {S}econd {C}onf., {U}niv. {O}klahoma, {N}orman, {O}kla., 1975)}, volume Vol. 26 of {\em Lect. Notes Pure Appl. Math.}, pages 171--223. Dekker, New York-Basel, 1977.

\bibitem[HRW99]{HRW99componentwiselinear}
J.~Herzog, V.~Reiner, and V.~Welker.
\newblock Componentwise linear ideals and {G}olod rings.
\newblock {\em Michigan Math. J.}, 46(2):211--223, 1999.

\bibitem[HVT10]{HT2010pathideals}
Jing He and Adam Van~Tuyl.
\newblock Algebraic properties of the path ideal of a tree.
\newblock {\em Comm. Algebra}, 38(5):1725--1742, 2010.

\bibitem[Jon08]{jonsson08book}
Jakob Jonsson.
\newblock {\em Simplicial complexes of graphs}, volume 1928 of {\em Lecture Notes in Mathematics}.
\newblock Springer-Verlag, Berlin, 2008.

\bibitem[JV25]{JV25dimension2posets}
Rizwan Jahangir and Dharm Veer.
\newblock On {C}ohen-{M}acaulay posets of dimension two and permutation graphs.
\newblock {\em Bull. Malays. Math. Sci. Soc.}, 48(4):Paper No. 133, 10, 2025.

\bibitem[JZ10]{JZ}
A.~Soleyman Jahan and X.~Zheng.
\newblock Ideals with linear quotients.
\newblock {\em Journal of Combinatorial Theory, Series A}, 117(1):104--110, 2010.

\bibitem[Kim22]{kim22cyclelength3}
Jinha Kim.
\newblock The homotopy type of the independence complex of graphs with no induced cycles of length divisible by 3.
\newblock {\em European J. Combin.}, 104:Paper No. 103534, 9, 2022.

\bibitem[KM06]{KR06lerayregularity}
Gil Kalai and Roy Meshulam.
\newblock Intersections of {L}eray complexes and regularity of monomial ideals.
\newblock {\em J. Combin. Theory Ser. A}, 113(7):1586--1592, 2006.

\bibitem[KNQ26]{KNQ24squarefree}
Elshani Kamberi, Francesco Navarra, and Ayesha~Asloob Qureshi.
\newblock On squarefree powers of simplicial trees.
\newblock {\em J. Algebra}, 693:240--276, 2026.

\bibitem[Koz07]{kozlov2007combinatorial}
D.~Kozlov.
\newblock {\em Combinatorial Algebraic Topology}.
\newblock Algorithms and Computation in Mathematics. Springer Berlin Heidelberg, 2007.

\bibitem[Mat19]{Matsushita19matching}
Takahiro Matsushita.
\newblock Matching complexes of small grids.
\newblock {\em Electron. J. Combin.}, 26(3):Paper No. 3.1, 8, 2019.

\bibitem[PB80]{PB}
J.~Scott Provan and Louis~J. Billera.
\newblock Decompositions of simplicial complexes related to diameters of convex polyhedra.
\newblock {\em Mathematics of Operations Research}, 5(4):576--594, 1980.

\bibitem[Rei76]{R}
Gerald~Allen Reisner.
\newblock Cohen-macaulay quotients of polynomial rings.
\newblock {\em Advances in Mathematics}, 21(1):30--49, 1976.

\bibitem[Sin20]{anurag20forest}
Anurag Singh.
\newblock Bounded degree complexes of forests.
\newblock {\em Discrete Math.}, 343(10):112009, 7, 2020.

\bibitem[Ter99]{terai99Alexanderduality}
Naoki Terai.
\newblock Alexander duality theorem and {S}tanley-{R}eisner rings.
\newblock Number 1078, pages 174--184. 1999.
\newblock Free resolutions of coordinate rings of projective varieties and related topics (Japanese) (Kyoto, 1998).

\bibitem[{The}20]{sage}
{The Sage Developers}.
\newblock {\em {S}ageMath, the {S}age {M}athematics {S}oftware {S}ystem ({V}ersion 9.2)}, 2020.
\newblock {\tt https://www.sagemath.org}.

\end{thebibliography}

\end{document}